\documentclass[10pt,reqno,draft,a4paper]{amsart}
\usepackage{tikz}
\usepackage{enumerate}
\usepackage{amsmath,amssymb,amsfonts} 
\usepackage{amsthm}
\usepackage{mathptmx}
\usepackage{helvet}
\usepackage{courier}
\usepackage{microtype}
\usepackage{hyperref}
\usepackage{nicefrac}
\usetikzlibrary{arrows,snakes}

\newtheorem{theorem}{Theorem}
\newtheorem{proposition}[theorem]{Proposition}
\newtheorem{lemma}[theorem]{Lemma}
\theoremstyle{definition}

\theoremstyle{remark}
\newtheorem{example}{Example}
\newcommand{\cptl}{{\mathcal{K}}} 
\newcommand{\strat}{{\mathcal{P}}}
\newcommand{\selec}{{\mathcal{S}}} 
\newcommand{\process}{{\mathcal{F}}}

\newcommand{\gamble}{{\mathcal{G}}} 
 
\newcommand{\smove}{{\mathbf{s}}} 
\newcommand{\wmove}{{\mathbf{w}}}
\newcommand{\smoves}{{\mathbf{S}}}
\newcommand{\wmoves}{{\mathbf{W}}}
\newcommand{\init}{\square} 
\newcommand{\fins}{\Omega}
\newcommand{\fin}{\omega} 
\newcommand{\sits}{\Omega^\lozenge}
\newcommand{\nonfins}{\sits\setminus\fins} 
\newcommand{\gain}{\lambda}
\newcommand{\price}{\mathbb{E}} 
\newcommand{\lprice}{\underline{\price}}
\newcommand{\uprice}{\overline{\price}}
\newcommand{\set}[2]{{\left\{#1\colon#2\right\}}}
\newcommand{\abs}[1]{{\lvert#1\rvert}}
\newcommand{\propset}[1]{{\left\{#1\right\}}} 
\newcommand{\pspace}{\Omega}
\newcommand{\posty}{\omega} 
\newcommand{\rdesirs}{\mathcal{R}}

\newcommand{\natex}{\mathcal{E}}
\newcommand{\gambles}{\mathcal{G}}  
\newcommand{\solp}{\mathcal{M}}
\newcommand{\partit}{\mathcal{B}}
\newcommand{\pr}{P} 
\newcommand{\lpr}{\underline{\pr}}
\newcommand{\upr}{\overline{\pr}}
 
\newcommand{\upset}[1]{{{\uparrow}#1}}
\newcommand{\precedes}{\sqsubseteq} 
\newcommand{\sprecedes}{\sqsubset}
\newcommand{\follows}{\sqsupseteq}
\newcommand{\sfollows}{\sqsupset}
\newcommand{\reals}{\mathbb{R}}
\newcommand{\heads}{h}
\newcommand{\tails}{t}
\newcommand{\flips}{\mathcal{F}}
\newcommand{\nheads}{\mathcal{N}}
\newcommand{\states}{\mathcal{X}}
\newcommand{\ntuple}[2]{{#1}_1,\dots,{#1}_{#2}}
\newcommand{\tuple}[2]{(\ntuple{#1}{#2})}
\DeclareMathOperator{\extremes}{ext}

\begin{document}
\title{Imprecise probability trees: Bridging two theories of imprecise
  probability} 
\author{Gert de Cooman \and Filip Hermans}
\address{Ghent University, SYSTeMS Research Group, Technologiepark --
  Zwijnaarde 914, 9052 Zwijnaarde, Belgium}
\email{\{gert.decooman,filip.hermans\}@UGent.be}

\begin{abstract}
  We give an overview of two approaches to probability theory where lower and upper   probabilities, rather than probabilities, are used: Walley's behavioural theory of imprecise   probabilities, and Shafer and Vovk's game-theoretic account of probability. We show that the   two theories are more closely related than would be suspected at first sight, and we establish   a correspondence between them that (i) has an interesting interpretation, and (ii) allows us to   freely import results from one theory into the other.  Our approach leads to an account of   probability trees and random processes in the framework of Walley's theory. We indicate how our   results can be used to reduce the computational complexity of dealing with imprecision in   probability trees, and we prove an interesting and quite general version of the weak law of   large numbers.
\end{abstract}

\maketitle

\section{Introduction}
In recent years, we have witnessed the growth of a number of theories of uncertainty, where imprecise (lower and upper) probabilities and previsions, rather than precise (or point-valued) probabilities and previsions, have a central part. Here we consider two of them, Glenn Shafer and Vladimir Vovk's game-theoretic account of probability \cite{shafer2001}, which is introduced in Section~\ref{sec:shafer-and-vovk}, and Peter Walley's behavioural theory \cite{walley1991}, outlined in Section~\ref{sec:walley}. These seem to have a rather different interpretation, and they certainly have been influenced by different schools of thought: Walley follows the tradition of Frank Ramsey \cite{ramsey1931}, Bruno de Finetti \cite{finetti19745} and Peter Williams \cite{williams2007} in trying to establish a rational model for a subject's beliefs in terms of her behaviour. Shafer and Vovk follow an approach that has many other influences as well, and is strongly coloured by ideas about gambling systems and martingales. They use Cournot's Principle to interpret lower and upper probabilities (see \cite{shafer2003}; and \cite[Chapter~2]{shafer2001} for a nice historical overview), whereas on Walley's approach, lower and upper probabilities are defined in terms of a subject's betting rates.
\par
What we set out to do here,\footnote{An earlier and condensed version of this paper, with much less   discussion and without proofs, was presented at the ISIPTA '07 conference \cite{cooman2007b}.} and in particular in Sections~\ref{sec:connections} and~\ref{sec:interpretation}, is to show that in many practical situations, the two approaches are strongly connected.\footnote{Our line of reasoning here   should be contrasted with the one in \cite{shafer2003}, where Shafer \textit{et al.} use the   game-theoretic framework developed in \cite{shafer2001} to construct a theory of predictive upper   and lower previsions whose interpretation is based on Cournot's Principle. See also the comments   near the end of Section~\ref{sec:interpretation}.}  This implies that quite a few results, valid in one theory, can automatically be converted and reinterpreted in terms of the other.  Moreover, we shall see that we can develop an account of coherent immediate prediction in the context of Walley's behavioural theory, and prove, in Section~\ref{sec:weak-law}, a weak law of large numbers with an intuitively appealing interpretation. We use this weak law in Section~\ref{sec:scoring} to suggest a way of scoring a predictive model that satisfies A.~Philip Dawid's \emph{Prequential Principle} \cite{dawid1984,dawid1999}.
\par
Why do we believe these results to be important, or even relevant, to AI? Probabilistic models are intended to represent an agent's beliefs about the world he is operating in, and which describe and even determine the actions he will take in a diversity of situations. Probability theory provides a normative system for reasoning and making decisions in the face of uncertainty. Bayesian, or precise, probability models have the property that they are completely decisive: a Bayesian agent always has an optimal choice when faced with a number of alternatives, whatever his state of information. While many may view this as an advantage, it is not always very realistic. Imprecise probability models try to deal with this problem by explicitly allowing for indecision, while retaining the normative, or coherentist stance of the Bayesian approach. We refer to \cite{cooman2005c,walley1991,walley1996} for discussions about how this can be done.
\par
Imprecise probability models appear in a number of AI-related fields. For instance in \emph{probabilistic logic}: it was already known to George Boole \cite{boole1847} that the result of probabilistic inferences may be a set of probabilities (an imprecise probability model), rather than a single probability. This is also important for dealing with missing or incomplete data, leading to so-called partial identification of probabilities, see for instance \cite{cooman2004b,manski2003}. There is also a growing literature on so-called \emph{credal nets} \cite{cozman2000,cozman2005}: these are essentially Bayesian nets with imprecise conditional probabilities. 
\par
We are convinced that it is mainly the mathematical and computational complexity often associated with imprecise probability models that is keeping them from becoming a more widely used tool for modelling uncertainty. But we believe that the results reported here can help make inroads in reducing this complexity. Indeed, the upshot of our being able to connect Walley's approach with Shafer and Vovk's, is twofold. First of all, we can develop a theory of \emph{imprecise probability trees}: probability trees where the transition from a node to its children is described by an imprecise probability model in Walley's sense. Our results provide the necessary apparatus for making inferences in such trees. And because probability trees are so closely related to random processes, this effectively brings us into a position to start developing a theory of (event-driven) random processes where the uncertainty can be described using imprecise probability models. We illustrate this in Examples~\ref{ex:coins} and~\ref{ex:many-coins}, and in Section~\ref{sec:backwards-recursion}. 
\par
Secondly, we are able to prove so-called Marginal Extension results (Theorems~\ref{theo:natex} and~\ref{theo:concatenation}, Proposition~\ref{prop:local-models}), which lead to backwards recursion, and dynamic programming-like methods that allow for an exponential reduction in the computational complexity of making inferences in such imprecise probability trees. This is also illustrated in Examples~\ref{ex:many-coins} and Section~\ref{sec:backwards-recursion}. For (precise) probability trees, similar techniques were described in Shafer's book on causal reasoning \cite{shafer1996}. They seem to go back to Christiaan Huygens, who drew the first probability tree, and showed how to reason with it, in his solution to Pascal and Fermat's Problem of Points.\footnote{See Section~\ref{sec:backwards-recursion} for more details and precise references.}

\section{Shafer and Vovk's game-theoretic approach to probability}
\label{sec:shafer-and-vovk}
In their game-theoretic approach to probability \cite{shafer2001}, Shafer and Vovk consider a game with two players, Reality and Sceptic, who play according to a certain \emph{protocol}. They obtain the most interesting results for what they call \emph{coherent   probability protocols}. This section is devoted to explaining what this means.

\subsection{Reality's event tree}
We begin with a first and basic assumption, dealing with how the first player, Reality, plays.
\begin{enumerate}[G1.]
\item Reality makes a number of moves, where the possible next moves may depend on the previous moves   he has made, but do not in any way depend on the previous moves made by Sceptic.
\end{enumerate}
This means that we can represent his game-play by an event tree (see also \cite{shafer1996a,shafer2000} for more information about event trees). We restrict ourselves here to the discussion of \emph{bounded protocols}, where Reality makes only a finite and bounded number of moves from the beginning to the end of the game, whatever happens.  But we don't exclude the possibility that at some point in the tree, Reality has the choice between an infinite number of next moves. We shall come back to these assumptions further on, once we have the appropriate notational tools to make them more explicit.\footnote{Essentially, the width of the tree may be infinite, but its depth should be finite.}
\par
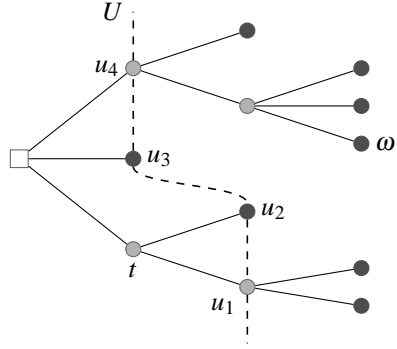
\begin{figure}[h]
  \centering 
  \tikzstyle{terminal}=[circle,draw=black!70,fill=black!70] 
  \tikzstyle{root}=[rectangle,draw=black!60]   
  \tikzstyle{nonterminal}=[circle,draw=black!60,fill=black!30,minimum size=.5mm] 
  \tikzstyle{level 1}=[sibling distance=12mm]   
  \tikzstyle{level 2}=[sibling distance=10mm] 
  \tikzstyle{level 3}=[sibling distance=5mm]
  \tikzstyle{cut}=[semithick,dashed,draw=black]
  \begin{tikzpicture}
    \node[root] (root) {} [grow=right,inner sep=.7mm] 
    child {node[nonterminal] [label=below:$t$] {} 
      child {node[nonterminal][label=below left:$u_1$] (cut1) {} 
        child {node[terminal] {}} 
        child {node[terminal] {}} } 
      child {node[terminal][label=right:$u_2$] (cut2) {}} } 
    child {node[terminal] [label=right:$u_3$] (cut3) {} } 
    child {node[nonterminal] [label=left:$u_4$] (cut4) {} 
      child {node[nonterminal]{} 
        child {node[terminal] [label=right:$\fin$] {}} 
        child {node[terminal] {}} 
        child {node[terminal] {}} } 
      child {node[terminal] {}}} ;
    \draw[cut] (cut3) -- (cut4) -- +(0,.75) node[left] {$U$};
    \draw[cut] (cut2) -- (cut1) -- +(0,-.75);
    \draw[cut] (cut3.south) .. controls +(down:.25) and +(up:.25) .. (cut2.north);
  \end{tikzpicture}
  \caption{A simple event tree for Reality, displaying the initial situation $\init$, other     non-terminal situations (such as $t$) as grey circles, and paths, or terminal situations,     (such as $\fin$) as black circles. Also depicted is a cut $U=\{u_1,u_2,u_3,u_4\}$ of $\init$.     Observe that $t$ (strictly) precedes $u_1$: $t\sprecedes u_1$, and that $C(t)=\{u_1,u_2\}$ is     the children cut of $t$.}
  \label{fig:tree-one}
\end{figure}
\par
Let us establish some terminology related to Reality's event tree. 

\subsubsection{Paths, situations and events}
A \emph{path} in the tree represents a possible sequence of moves for Reality from the beginning to the end of the game. We denote the set of all possible paths $\fin$ by $\fins$, the \emph{sample space} of the game.  
\par
A \emph{situation} $t$ is some connected segment of a path that is \emph{initial}, i.e., starts at the root of the tree.  It identifies the moves Reality has made up to a certain point, and it can be identified with a node in the tree. We denote the set of all situations by $\sits$. It includes the set $\fins$ of \emph{terminal} situations, which can be identified with paths. All other situations are called \emph{non-terminal}; among them is the \emph{initial} situation $\init$, which represents the empty initial segment. See Fig.~\ref{fig:tree-one} for a simple graphical example explaining these notions.
\par
If for two situations $s$ and $t$, $s$ is a(n initial) segment of $t$, then we say that $s$ \emph{precedes} $t$ or that $t$ \emph{follows} $s$, and write $s\precedes t$, or alternatively $t\follows s$. If $\fin$ is a path and $t\precedes\fin$ then we say that the path $\fin$ \emph{goes through} situation $t$. We write $s\sprecedes t$, and say that \emph{$s$ strictly precedes $t$}, if $s\precedes t$ and $s\neq t$.
\par
An \emph{event} $A$ is a set of paths, or in other words, a subset of the sample space: $A\subseteq\fins$. With an event $A$, we can associate its \emph{indicator} $I_A$, which is the real-valued map on $\fins$ that assumes the value $1$ on $A$, and $0$ elsewhere.
\par 
We denote by $\upset{t}:=\set{\fin\in\fins}{t\precedes\fin}$ the set of all paths that go through $t$: $\upset{t}$ is the event that corresponds to Reality getting to a situation $t$.  It is clear that not all events will be of the type $\upset{t}$. Shafer \cite{shafer1996} calls events of this type \emph{exact}. Further on, in Section~\ref{sec:connections}, exact events will be the only events that can be legitimately conditioned on, because they are the only events that can be foreseen may occur as part of Reality's game-play.

\subsubsection{Cuts of a situation}
Call a \emph{cut} $U$ of a situation $t$ any set of situations that 
follow $t$, and such that for all paths $\fin$ through $t$, there is a unique $u\in U$ that $\fin$ goes through. In other words:
\begin{enumerate}[(i)]
\item $(\forall u\in U)(u\follows t)$; and
\item $(\forall\fin\follows t)(\exists!u\in U)(\fin\follows u)$; 
\end{enumerate}
see also Fig.~\ref{fig:tree-one}. Alternatively, a set $U$ of situations is a cut of $t$ if and only if the corresponding set $\set{\upset{u}}{u\in U}$ of exact events is a partition of the exact event $\upset{t}$. A cut can be interpreted as a (complete) stopping time.
\par
If a situation $s\follows t$ precedes (follows) some element of a cut $U$ of $t$, then we say that $s$ \emph{precedes} (\emph{follows}) $U$, and we write $s\precedes U$ ($s\follows U$). Similarly for `strictly precedes (follows)'.  For two cuts $U$ and $V$ of $t$, we say that $U$ \emph{precedes} $V$ if each element of $U$ is followed by some element of $V$.
\par
A \emph{child} of a non-terminal situation $t$ is a situation that immediately follows it. The set $C(t)$ of children of $t$ constitutes a cut of $t$, called its \emph{children cut}. Also, the set $\fins$ of terminal situations is a cut of $\init$, called its \emph{terminal cut}. The event $\upset{t}$ is the corresponding terminal cut of a situation $t$.

\subsubsection{Reality's move spaces}
\begin{figure}[h]
  \centering
  \begin{tikzpicture}[level distance=30mm]
    \tikzstyle{terminal}=[circle,draw=black!70,fill=black!70] 
    \tikzstyle{root}=[rectangle,draw=black!60]   
    \tikzstyle{nonterminal}=[circle,draw=black!60,fill=black!30,minimum size=.5mm] 
    \tikzstyle{level 1}=[sibling distance=15mm] 
    \tikzstyle{level 2}=[sibling distance=15mm] 
    \tikzstyle{level 3}=[sibling distance=15mm]
    \node[root] (root) {} [grow=right,inner sep=.7mm] 
    child {node[nonterminal,label=below:$t$] (t) {}
      child {node[nonterminal,label=above right:$t\wmove_2$] (tw2) {}
        child {node[terminal] {}}
        child {node[terminal] {}}
        edge from parent {node[below] {$\wmove_2$}}}
      child {node[terminal,label=above right:$t\wmove_1$] (tw1) {}
        edge from parent {node[above] {$\wmove_1$}}}}
    child {node[terminal] {} };
    \draw[semithick,dashed] (tw1) -- +(0,.85);
    \draw[semithick,dashed] (tw1) -- (tw2) -- +(0,-1) node[right] {$C(t)$};
    \draw[darkgray,dashed] (t) -- +(-2,-1) node[below] {$\wmoves_t=\{\wmove_1,\wmove_2\}$}; 
  \end{tikzpicture}
  \caption{An event tree for Reality, with the move space $\wmoves_t$ and the corresponding children cut $C(t)$ of a non-terminal situation $t$.}
  \label{fig:moves} 
\end{figure}
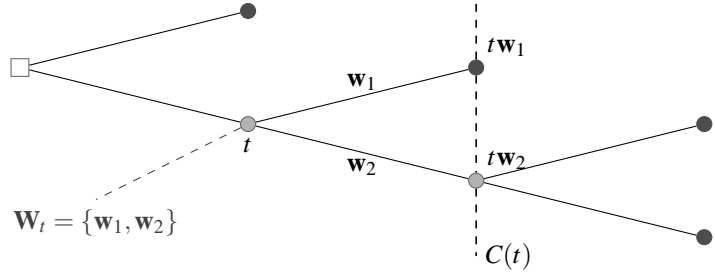

We call a \emph{move} $\wmove$ for Reality in a non-terminal situation $t$ an arc that connects $t$ with one of its children $s\in C(t)$, meaning that $s=t\wmove$ is the concatenation of the segment $t$ and the arc $\wmove$. See Fig.~\ref{fig:moves}.  
\par
\par
Reality's \emph{move space} in $t$ is the set $\wmoves_t$ of those moves $\wmove$ that Reality can make in $t$: $\wmoves_t=\set{\wmove}{t\wmove\in C(t)}$. We have already mentioned that $\wmoves_t$ may be (countably or uncountably) infinite: there may be situations where reality has the choice between an infinity of next moves. But every $\wmoves_t$ should contain at least two elements:  otherwise there is no choice for Reality to make in situation $t$.

\subsection{Processes and variables}
We now have all the necessary tools to represent Reality's game-play. This game-play can be seen as a basis for an \emph{event-driven}, rather than a time-driven, account of a theory of uncertain, or random, processes. The driving events are, of course, the moves that Reality makes.\footnote{These so-called \emph{Humean} events shouldn't be confused with the \emph{Moivrean} events we have considered before, and which are subsets of the sample space $\fins$. See Shafer \cite[Chapter~1]{shafer1996} for terminology and more explanation.} In a theory of processes, we generally consider things that depend on (the succession of) these moves. This leads to the following definitions. 
\par
Any (partial) function on the set of situations $\sits$ is called a \emph{process}, and any process whose domain includes all situations that follow a situation $t$ is called a \emph{$t$-process}. Of course, a $t$-process is also an $s$-process for all $s\follows t$; when we call it an $s$-process, this means that we are restricting our attention to its values in all situations that follow $s$.
\par
A special example of a  $t$-process is the \emph{distance} $d(t,\cdot)$ which for any situation $s\follows t$ returns the number of steps $d(t,s)$ along the tree from $t$ to $s$. When we said before that we are only considering \emph{bounded protocols}, we meant that there is a natural number $D$ such that $d(t,s)\leq D$ for all situations $t$ and all $s\follows t$.
\par
Similarly, any (partial) function on the set of paths $\fins$ is called a \emph{variable}, and any variable on $\fins$ whose domain includes all paths that go through a situation $t$ is called a \emph{$t$-variable}. If we restrict a $t$-process $\process$ to the set $\upset{t}$ of all terminal situations that follow $t$, we obtain a $t$-variable, which we denote by $\process_\fins$.
\par
If $U$ is a cut of $t$, then we call a $t$-variable $g$ \emph{$U$-measurable} if for all $u$ in $U$, $g$ assumes the same value $g(u):=g(\fin)$ for all paths $\fin$ that go through $u$. In that case we can also consider $g$ as a variable on $U$, which we denote as $g_U$.
\par
If $\process$ is a $t$-process, then with any cut $U$ of $t$ we can associate a $t$-variable $\process_U$, which assumes the same value $\process_U(\fin):=\process(u)$ in all $\fin$ that follow $u\in U$. This $t$-variable is clearly $U$-measurable, and can be considered as a variable on $U$. This notation is consistent with the notation $\process_\fins$ introduced earlier.
\par
Similarly, we can associate with $\process$ a new, \emph{$U$-stopped}, $t$-process $U(\process)$, as follows:
\begin{equation*}
  U(\process)(s):=
  \begin{cases}
    \process(s)&\text{if $t\precedes s\precedes U$}\\
    \process(u)&\text{if $u\in U$ and $u\precedes s$}.
  \end{cases}
\end{equation*}
The $t$-variable $U(\process)_\fins$ is $U$-measurable, and is actually equal to $\process_U$:
\begin{equation}\label{eq:stopped}
  U(\process)_\fins=\process_U.
\end{equation}
The following intuitive example will clarify these notions.

\begin{example}[Flipping coins]\label{ex:coins}
  Consider flipping two coins, one after the other. This leads to the event tree depicted in Fig.~\ref{fig:coins}. The identifying labels for the situations should be intuitively clear: e.g., in the initial situation `$\init=?,?$' none of the coins have been flipped, in the non-terminal situation `$\heads,?$' the first coin has landed `heads' and the second coin hasn't been flipped yet, and in the terminal situation `$\tails,\tails$' both coins have been flipped and have landed `tails'.
\par
\begin{figure}[h]
  \centering
  \begin{tikzpicture}[level distance=30mm]
    \tikzstyle{first}=[rectangle,rounded corners,draw=black!50,fill=black!20]
    \tikzstyle{root}=[rectangle,rounded corners,draw=black!50,fill=black!20]
    \tikzstyle{second}=[rectangle,rounded corners,draw=black!50,fill=black!20]
    \tikzstyle{level 1}=[sibling distance=20mm]
    \tikzstyle{level 2}=[sibling distance=10mm]
    \tikzstyle{cut1}=[semithick,dashed,draw=black]
    \tikzstyle{cut2}=[semithick,dotted,draw=black]
    \node[root] (root) {$?,?$} [grow=right,inner sep=.7mm]
    child {node[first] (t) {$\tails,?$}
      child {node[second] (tt) {$\tails,\tails$}}
      child {node[second] (th) {$\tails,\heads$}}}
    child {node[first] (h) {$\heads,?$}
      child {node[second] (ht) {$\heads,\tails$}}
      child {node[second] (hh) {$\heads,\heads$}}};
    \node (in) [above left of=h] {};
    \node (uit) [below right of=tt] {};
    \draw[cut1] (h) -- +(0,1);
    \draw[cut1] (h) -- (t) -- +(0,-1.5) node[left] {$X^1$};
    \draw[cut2] (in.east) to [bend left] (h.north);
    \draw[cut2] (h.south) .. controls +(down:.5) and +(up:.5) .. (th.north);
    \draw[cut2] (th) -- (tt) -- +(0,-1) node[left] {$U$};
  \end{tikzpicture}
  \caption{The event tree associated with two successive coin flips. Also depicted are two cuts, $X_1$ and $U$, of the initial situation.}
  \label{fig:coins}
\end{figure}
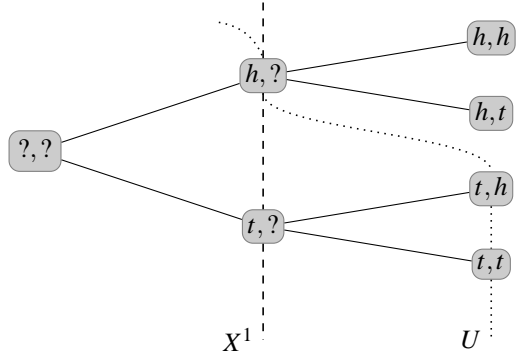
\par
First, consider the real process $\nheads$, which in each situation $s$, returns the number $\nheads(s)$ of heads obtained so far, e.g., $\nheads(?,?)=0$ and $\nheads(\heads,?)=1$. If we restrict the process $\nheads$ to the set $\fins$ of all terminal elements, we get a real variable $\nheads_\fins$, whose values are: $\nheads_\fins(\heads,\heads)=2$, $\nheads_\fins(\heads,\tails)=\nheads_\fins(\tails,\heads)=1$ and $\nheads_\fins(\tails,\tails)=0$.
\par
Consider the cut $U$ of the initial situation, which corresponds to the following stopping time: ``stop after two flips, or as soon as an outcome is heads''; see Fig.~\ref{fig:coins}. The values of the corresponding variable $\nheads_U$ are given by: $\nheads_U(\heads,\heads)=\nheads_U(\heads,\tails)=1$, $\nheads_U(\tails,\heads)=1$ and $\nheads_U(\tails,\tails)=0$. So $\nheads_U$ is $U$-measurable, and can therefore be considered as a map on the elements $\heads,?$ and $\tails,\heads$ and $\tails,\tails$ of $U$, with in particular $\nheads_U(\heads,?)=1$.
\par
Next, consider the processes $\flips,\flips^1,\flips^2\colon\sits\to\{\heads,\tails,?\}$, defined as follows:
\begin{center}
  \begin{tabular}[h]{c|ccccccc}
    $s$ & $?,?$ & $\heads,?$ & $\tails,?$ & $\heads,\heads$ 
    & $\heads,\tails$ & $\tails,\heads$ & $\tails,\tails$ \\
    $\flips(s)$ & $?$ & $\heads$ & $\tails$ & $\heads$ & $\tails$ & $\heads$ & $\tails$ \\
    $\flips^1(s)$ & $?$ & $\heads$ & $\tails$ & $\heads$ & $\heads$ & $\tails$ & $\tails$ \\
    $\flips^2(s)$ & $?$ & $?$ & $?$ & $\heads$ & $\tails$ & $\heads$ & $\tails$ 
  \end{tabular}
\end{center}
$\flips$ returns the outcome of the latest, $\flips^1$ the outcome of the first, and $\flips^2$ that of the second coin flip. The associated variables $\flips^1_\fins$ and $\flips^2_\fins$ give, in each element of the sample space, the respective outcomes of the first and second coin flips.
\par
The variable $\flips^1_\fins$ is $X^1$-measurable: as soon as we reach (any situation on) the cut $X_1$, its value is completely determined, i.e., we know the outcome of the first coin flip; see Fig.~\ref{fig:coins} for the definition of $X^1$.
\par
We can associate with the process $\flips$ the variable $\flips_{X^1}$ that is also $X^1$-measurable: it returns, in any element of the sample space, the outcome of the first coin flip. Alternatively, we can stop the process $\flips$ after one coin flip, which leads to the $X^1$-stopped process $X^1(\flips)$. This new process is of course equal to $\flips^1$, and for the corresponding variable $\flips^1_\fins$, we have that $X^1(\flips)_\fins=\flips^1_\fins=\flips_{X^1}$; also see Eq.~\eqref{eq:stopped}. $\blacklozenge$
\end{example}

\subsection{Sceptic's game-play}
We now turn to the other player, Sceptic. His possible moves may well depend on the previous moves that Reality has made, in the following sense. In each non-terminal situation $t$, he has some set $\smoves_t$ of moves $\smove$ available to him, called Sceptic's \emph{move space} in $t$. We make the following assumption:
\begin{enumerate}[G2.]
\item In each non-terminal situation $t$, there is a (positive or negative) gain for Sceptic associated with each of the possible   moves $\smove$ in $\smoves_t$ that Sceptic can make. This gain depends only on the situation $t$ and the next move $\wmove$ that   Reality will make.
\end{enumerate}
This means that for each non-terminal situation $t$ there is a \emph{gain function} $\gain_t\colon\smoves_t\times\wmoves_t\to\reals$, such that $\gain_t(\smove,\wmove)$ represents the change in Sceptic's capital in situation $t$ when he makes move $\smove$ and Reality makes move $\wmove$.

\subsubsection{Strategies and capital processes}
Let us introduce some further notions and terminology related to Sceptic's game-play. A \emph{strategy} $\strat$ for Sceptic is a partial process defined on the set $\nonfins$ of non-terminal situations, such that $\strat(t)\in\smoves_t$ is the corresponding move that Sceptic will make in each non-terminal situation $t$.
\par
With each such strategy $\strat$ there corresponds a \emph{capital process} $\cptl^\strat$, whose value in each situation $t$ gives us Sceptic's capital accumulated so far, when he starts out with zero capital in $\init$ and plays according to the strategy $\strat$. It is given by the recursion relation
\begin{equation*}
  \cptl^\strat(t\wmove)
  =\cptl^\strat(t)+\gain_t(\strat(t),\wmove),
  \quad\wmove\in\wmoves_t,
\end{equation*}
with initial condition $\cptl^\strat(\init)=0$. Of course, when Sceptic starts out (in $\init$) with capital $\alpha$ and uses strategy $\strat$, his corresponding accumulated capital is given by the process $\alpha+\cptl^\strat$. In the terminal situations, his accumulated capital is then given by the real variable $\alpha+\cptl^\strat_\fins$.
\par
If we start in a non-terminal situation $t$, rather than in $\init$, then we can consider $t$-strategies $\strat$ that tell Sceptic how to move starting from $t$ onwards, and the corresponding capital process $\cptl^\strat$ is then also a $t$-process, that tells us how much capital Sceptic has accumulated since starting with zero capital in situation $t$ and using $t$-strategy $\strat$.

\subsubsection{Lower and upper prices}
The assumptions G1 and G2 outlined above determine so-called \emph{gambling protocols}. They are sufficient for us to be able to define lower and upper prices for real variables.  
\par
Consider a non-terminal situation $t$ and a real $t$-variable $f$.  The \emph{upper price $\uprice_t(f)$ for $f$ in $t$} is defined as the infimum capital $\alpha$ that Sceptic has to start out with in $t$ in order that there would be some $t$-strategy $\strat$ such that his accumulated capital $\alpha+\cptl^\strat$ allows him, at the end of the game, to hedge $f$, whatever moves Reality makes after $t$:
\begin{equation}\label{eq:uprice}
  \uprice_t(f)
  :=\inf\set{\alpha}
  {\text{$\alpha+\cptl^\strat_\fins\geq f$ for some $t$-strategy $\strat$}},
\end{equation}
where $\alpha+\cptl^\strat_\fins\geq f$ is taken to mean that $\alpha+\cptl^\strat(\fin)\geq f(\fin)$ for all terminal situations $\fin$ that go through $t$. Similarly, for the \emph{lower price $\lprice_t(f)$ for $f$ in $t$}:
\begin{equation}\label{eq:lprice}
  \lprice_t(f)
  :=\sup\set{\alpha}
  {\text{$\alpha-\cptl^\strat_\fins\leq f$ for some $t$-strategy $\strat$}},
\end{equation}
so $\lprice_t(f)=-\uprice_t(-f)$. If we start from the initial situation $t=\square$, we simply get the \emph{upper and lower prices} for a real variable $f$, which we also denote by $\uprice(f)$ and $\lprice(f)$.

\subsubsection{Coherent probability protocols}
Requirements~G1 and~G2 for gambling protocols allow the moves, move spaces and gain functions for Sceptic to be just about anything. We now impose further conditions on Sceptic's move spaces.
\par 
A gambling protocol is called a \emph{probability protocol} when besides G1 and G2, two more requirements are satisfied.
\begin{enumerate}[P1.]
\item For each non-terminal situation $t$, Sceptic's move space $\smoves_t$ is a convex cone in some linear space:   $a_1\smove_1+a_2\smove_2\in\smoves_t$ for all non-negative real numbers $a_1$ and $a_2$ and all $\smove_1$ and   $\smove_2$ in $\smoves_t$.
\item For each non-terminal situation $t$, Sceptic's gain function $\gain_t$ has the following linearity   property: $\gain_t(a_1\smove_1+a_2\smove_2,\wmove) =a_1\gain_t(\smove_1,\wmove)+a_2\gain_t(\smove_2,\wmove)$   for all non-negative real numbers $a_1$ and $a_2$, all $\smove_1$ and $\smove_2$ in $\smoves_t$ and all   $\wmove$ in $\wmoves_t$.
\end{enumerate}
Finally, a probability protocol is called \emph{coherent}\footnote{For a discussion of the use of `coherent' here, we refer to   \cite[Appendix~C]{shafer2003}.} when moreover:
\begin{enumerate}[C.]
\item For each non-terminal situation $t$, and for each $\smove$ in $\smoves_t$ there is some $\wmove$ in   $\wmoves_t$ such that $\gain_t(\smove,\wmove)\leq0$.
\end{enumerate}
It is clear what this last requirement means: in each non-terminal situation, Reality has a strategy for playing from $t$ onwards such that Sceptic can't (strictly) increase his capital from $t$ onwards, whatever $t$-strategy he might use.
\par
For such coherent probability protocols, Shafer and Vovk prove a number of interesting properties for the corresponding lower (and upper) prices. We list a number of them here. For any real $t$-variable $f$, we can associate with a cut $U$ of $t$ another special $U$-measurable $t$-variable $\lprice_U$ by $\lprice_U(f)(\fin)=\lprice_u(f)$, for all paths $\fin$ through $t$, where $u$ is the unique situation in $U$ that $\fin$ goes through. For any two real $t$-variables $f_1$ and $f_2$, $f_1\leq f_2$ is taken to mean that $f_1(\fin)\leq f_2(\fin)$ for all paths $\omega$ that go through $t$.

\begin{proposition}[Properties of lower and upper prices in a coherent
  probability protocol
  \protect{\cite{shafer2001}}]\label{prop:shafer-and-vovk}
  Consider a coherent probability protocol, let\/ $t$ be a non-terminal situation, $f$, $f_1$ and\/ $f_2$ real\/ $t$-variables, and   $U$ a cut of\/ $t$.  Then
  \begin{enumerate}[1.]
  \item $\inf_{\fin\in\upset{t}}f(\fin)\leq\lprice_t(f) \leq\uprice_t(f)\leq\sup_{\fin\in\upset{t}}f(\fin)$ [convexity];
  \item $\lprice_t(f_1+f_2)\geq\lprice_t(f_1)+\lprice_t(f_2)$ [super-additivity];
  \item $\lprice_t(\lambda f)=\lambda\lprice_t(f)$ for all real $\lambda\geq0$ [non-negative homogeneity];
  \item $\lprice_t(f+\alpha)=\lprice_t(f)+\alpha$ for all real $\alpha$ [constant additivity];
  \item $\lprice_t(\alpha)=\alpha$ for all real $\alpha$ [normalisation];
  \item $f_1\leq f_2$ implies that\/ $\lprice_t(f_1)\leq\lprice_t(f_2)$ [monotonicity];
  \item $\lprice_t(f)=\lprice_t(\lprice_U(f))$ [law of iterated expectation].
  \end{enumerate}
\end{proposition}
What is more, Shafer and Vovk use specific instances of such coherent probability protocols to prove various limit theorems (such as the law of large numbers, the central limit theorem, the law of the iterated logarithm), from which they can derive, as special cases, the well-known measure-theoretic versions. We shall come back to this in Section~\ref{sec:weak-law}.
\par
The game-theoretic account of probability we have described so far, is very general. But it seems to pay little or no attention to \emph{beliefs} that Sceptic, or other, perhaps additional players in these games might entertain about how Reality will move through its event tree. This might seem strange, because at least according to the personalist and epistemicist school, probability is all about beliefs. In order to find out how we can incorporate beliefs into the game-theoretic framework, we now turn to Walley's imprecise probability models.

\section{Walley's behavioural approach to probability}
\label{sec:walley}
In his book on the behavioural theory of imprecise probabilities \cite{walley1991}, Walley considers many different types of related uncertainty models. We shall restrict ourselves here to the most general and most powerful one, which also turns out to be the easiest to explain, namely coherent sets of really desirable gambles; see also \cite{walley2000}.
\par
Consider a non-empty set $\pspace$ of possible alternatives $\posty$, only one of which actually obtains (or will obtain); we assume that it is possible, at least in principle, to determine which alternative does so. Also consider a subject who is uncertain about which possible alternative actually obtains (or will obtain). A \emph{gamble} on $\pspace$ is a real-valued map on $\pspace$, and it is interpreted as an uncertain reward, expressed in units of some predetermined linear utility scale: if $\omega$ actually obtains, then the reward is $f(\omega)$, which may be positive or negative. We use the notation $\gambles(\pspace)$ for the set of all gambles on $\pspace$. Walley \cite{walley1991} assumes gambles to be bounded. We make no such boundedness assumption here.\footnote{The concept of a really desirable gamble (at least formally) allows for such a generalisation,   because the coherence axioms for real desirability nowhere hinge on such a boundedness assumption, at least not   from a technical mathematical point of view.}
\par
If a subject \emph{accepts} a gamble $f$, this is taken to mean that she is willing to engage in the transaction where, (i) first it is determined which $\omega$ obtains, and (ii) then she receives the reward $f(\omega)$. We can try and model the subject's beliefs about $\pspace$ by considering which gambles she accepts.

\subsection{Coherent sets of really desirable gambles}
Suppose our subject specifies some set $\rdesirs$ of gambles she accepts, called a \emph{set of really desirable gambles}. Such a set is called \emph{coherent} if it satisfies the following \emph{rationality requirements}:
\begin{enumerate}[D1.]
\item if\/ $f<0$ then $f\not\in\rdesirs$ [avoiding partial loss];
\item if\/ $f\geq0$ then $f\in\rdesirs$ [accepting partial gain];
\item if\/ $f_1$ and\/ $f_2$ belong to $\rdesirs$ then their (point-wise) sum $f_1+f_2$ also belongs to $\rdesirs$   [combination];
\item if $f$\/ belongs to $\rdesirs$ then its (point-wise) scalar product\/ $\lambda f$ also belongs to $\rdesirs$ for all   non-negative real numbers $\lambda$ [scaling].
\end{enumerate}
Here `$f<0$' means `$f\leq0$ and not $f=0$'.  Walley has also argued that, besides D1--D4, sets of really desirable gambles should satisfy an additional axiom:
\begin{enumerate}[D5.]
\item $\rdesirs$ is $\partit$-conglomerable for any partition $\partit$ of $\pspace$: if $I_Bf\in\rdesirs$ for all   $B\in\partit$, then also $f\in\rdesirs$ [full conglomerability].
\end{enumerate}
When the set $\pspace$ is finite, all its partitions are finite too, and therefore full conglomerability becomes a direct consequence of the finitary combination axiom D3. But when $\pspace$ is infinite, its partitions may be infinite too, and then full conglomerability is a very strong additional requirement, that is not without controversy. If a model $\rdesirs$ is $\partit$-conglomerable, this means that certain inconsistency problems when conditioning on elements $B$ of $\partit$ are avoided; see \cite{walley1991} for more details and examples. Conglomerability of belief models wasn't required by forerunners of Walley, such as Williams \cite{williams2007},\footnote{Axioms related to (D1)--(D4), but not (D5), were actually suggested by Williams for   bounded gambles. But it seems that we need at least some weaker form of (D5), namely the cut conglomerability (D5')   considered further on, to derive our main results: Theorems~\ref{theo:natex} and~\ref{theo:matching}.} or de Finetti \cite{finetti19745}. While we agree with Walley that conglomerability is a desirable property for sets of really desirable gambles, we do not believe that \emph{full} conglomerability is always necessary: it seems that we only need to require conglomerability with respect to those partitions that we actually intend to condition our model on.\footnote{The view expressed here seems related to Shafer's, as sketched near the end of   \cite[Appendix~1]{shafer1985}.}  This is the path we shall follow in Section~\ref{sec:connections}.

\subsection{Conditional lower and upper previsions}\label{sec:lower-upper}
Given a coherent set of really desirable gambles, we can define conditional lower and upper previsions as follows: for any gamble $f$ and any non-empty subset $B$ of $\pspace$, with indicator $I_B$,
\begin{gather}
  \upr(f\vert B):=\inf\set{\alpha}{I_B(\alpha-f)\in\rdesirs}\label{eq:updef}\\
  \lpr(f\vert B):=\sup\set{\alpha}{I_B(f-\alpha)\in\rdesirs},\label{eq:lowdef}
\end{gather}
so $\lpr(f\vert B)=-\upr(-f\vert B)$, and \emph{the lower prevision $\lpr(f\vert B)$ of $f$, conditional on $B$} is the supremum price $\alpha$ for which the subject will buy the gamble $f$, i.e., accept the gamble $f-\alpha$, contingent on the occurrence of $B$. Similarly, \emph{the upper prevision $\upr(f\vert B)$ of $f$, conditional on $B$} is the infimum price $\alpha$ for which the subject will sell the gamble $f$, i.e., accept the gamble $\alpha-f$, contingent on the occurrence of $B$. 
\par
For any event $A$, we define the conditional lower probability $\lpr(A\vert B):=\lpr(I_A\vert B)$, i.e., the subject's supremum rate for betting on the event $A$, contingent on the occurrence of $B$, and similarly for $\upr(A\vert B):=\upr(I_A\vert B)$.
\par
We want to stress here that by its definition [Eq.~\eqref{eq:lowdef}], $\lpr(f\vert B)$ is a conditional lower prevision on what Walley \cite[Section~6.1]{walley1991} has called the \emph{contingent interpretation}: it is a supremum acceptable price for buying the gamble $f$ \emph{contingent} on the occurrence of $B$, meaning that the subject accepts the contingent gambles $I_B(f-\lpr(f\vert B)+\epsilon)$, $\epsilon>0$, which are called off unless $B$ occurs. This should be contrasted with the \emph{updating interpretation} for the conditional lower prevision $\lpr(f\vert B)$, which is a subject's \emph{present} (before the occurrence of $B$) supremum acceptable price for buying $f$ after receiving the information that $B$ has occurred (and nothing else!). Walley's \emph{Updating Principle} \cite[Section~6.1.6]{walley1991}, which we shall accept, and use further on in Section~\ref{sec:connections}, (essentially) states that conditional lower previsions should be the same on both interpretations. There is also a third way of looking at a conditional lower prevision $\lpr(f\vert B)$, which we shall call the \emph{dynamic interpretation}, and where $\lpr(f\vert B)$ stands for the subject's supremum acceptable buying price for $f$ \emph{after she gets to know} $B$ has occurred. For precise conditional previsions, this last interpretation seems to be the one considered in \cite{goldstein1983,shafer1982,shafer1983,shafer2003}. It is far from obvious that there should be a relation between the first two and the third interpretations.\footnote{In \cite{shafer2003}, the authors seem to confuse the updating interpretation with the dynamic   interpretation when they claim that ``[their new understanding of lower and upper previsions]   justifies Peter Walley's updating principle''.}  We shall briefly come back to this distinction in the following sections.
\par
For any partition $\partit$ of $\pspace$, we let $\lpr(f\vert\partit):=\sum_{B\in\partit}I_B\lpr(f\vert B)$ be the gamble on $\pspace$ that in any element $\posty$ of $B$ assumes the value $\lpr(f\vert B)$, where $B$ is any element of $\partit$.
\par
The following properties of conditional lower and upper previsions associated with a coherent set of really desirable bounded gambles were (essentially) proved by Walley \cite{walley1991}, and by Williams \cite{williams2007}. We give the extension to potentially unbounded gambles:

\begin{proposition}[Properties of conditional lower and upper previsions
  \protect{\cite{walley1991}}]\label{prop:walley}
  Consider a coherent set of really desirable gambles $\rdesirs$, let\/ $B$ be any non-empty   subset of\/ $\pspace$, and let\/ $f$, $f_1$ and $f_2$ be gambles on $\pspace$.   Then\footnote{Here, as in Proposition~\ref{prop:shafer-and-vovk}, we implicitly assume that     whatever we write down is well-defined, meaning that for instance no sums of $-\infty$ and     $+\infty$ appear, and that the function $\lpr(f\vert\partit)$ is real-valued, and nowhere     infinite.  Shafer and Vovk don't seem to mention the need for this.\label{fn:well-defined}}
  \begin{enumerate}[1.]
  \item $\inf_{\posty\in B}f(\posty)\leq\lpr(f\vert B)\leq\upr(f\vert B)\leq\sup_{\posty\in B}f(\posty)$ [convexity];
  \item $\lpr(f_1+f_2\vert B)\geq\lpr(f_1\vert B)+\lpr(f_2\vert B)$ [super-additivity];
  \item $\lpr(\lambda f\vert B)=\lambda\lpr(f\vert B)$ for all real $\lambda\geq0$ [non-negative homogeneity];
  \item $\lpr(f+\alpha\vert B)=\lpr(f\vert B)+\alpha$ for all real $\alpha$ [constant additivity];
  \item $\lpr(\alpha\vert B)=\alpha$ for all real $\alpha$ [normalisation];
  \item $f_1\leq f_2$ implies that\/ $\lpr(f_1\vert B)\leq\lpr(f_2\vert B)$ [monotonicity];
  \item if\/ $\partit$ is any partition of\/ $\pspace$ that refines the partition $\{B,B^c\}$ and\/ $\rdesirs$ is $\partit$-conglomerable, then $\lpr(f\vert B)\geq\lpr(\lpr(f\vert\partit)\vert B)$ [conglomerative property].
  \end{enumerate}
\end{proposition}

The analogy between Propositions~\ref{prop:shafer-and-vovk} and~\ref{prop:walley} is striking, even if there is an equality in Proposition~\ref{prop:shafer-and-vovk}.7, where we have only an inequality in Proposition~\ref{prop:walley}.7.\footnote{Concatenation inequalities for lower prices do appear in the more general   context described in \cite{shafer2003}.} In the next section, we set out to identify the exact correspondence between the two models. We shall find a specific situation where applying Walley's theory leads to equalities rather than the more general inequalities of Proposition~\ref{prop:walley}.7.\footnote{This seems to   happen generally for what is called \emph{marginal extension} in a situation of immediate prediction, meaning that   we start out with, and extend, an initial model where we condition on increasingly finer partitions, and where the   initial conditional model for any partition deals with gambles that are measurable with respect to the finer   partitions; see \cite[Theorem~6.7.2]{walley1991} and \cite{miranda2006b}.}
\par
We now show that there can indeed be a strict inequality in Proposition~\ref{prop:walley}.7.

\begin{example}
  Consider an urn with red, green and blue balls, from which a ball will be drawn at random. Our subject   is uncertain about the colour of this ball, so $\pspace=\{r,g,b\}$. Assume that she assesses that she is willing to bet on this colour being red at rates up to (and including) $\nicefrac{1}{4}$, i.e., that she accepts the gamble $I_{\{r\}}-\nicefrac{1}{4}$. Similarly for the other two colours, so she also accepts the gambles $I_{\{g\}}-\nicefrac{1}{4}$ and $I_{\{b\}}-\nicefrac{1}{4}$. It is not difficult to prove using the coherence requirements D1--D4 and Eq.~\eqref{eq:lowdef} that the smallest coherent set of really desirable gambles $\rdesirs$ that includes these assessments satisfies $f\in\rdesirs\Leftrightarrow\lpr(f)\geq0$, where
\begin{equation*}
  \lpr(f)=\frac{3}{4}\frac{f(r)+f(g)+f(b)}{3}+\frac{1}{4}\min\{f(r),f(g),f(b)\}.  
\end{equation*}
For the partition $\partit=\{b,\{r,g\}\}$ (a Daltonist has observed the colour of the ball and tells the subject about it), it follows from Eq.~\eqref{eq:lowdef} after some manipulations that
\begin{equation*}
  \lpr(f\vert\{b\})=f(b)
  \text{ and }
  \lpr(f\vert\{r,g\})=\frac{2}{3}\frac{f(r)+f(g)}{2}+\frac{1}{3}\min\{f(r),f(g)\}.
\end{equation*}
If we consider $f=I_{\{g\}}$, then in particular $\lpr(\{g\}\vert\{b\})=0$ and $\lpr(\{g\}\vert\{r,g\})=\nicefrac{1}{3}$, so $\lpr(\{g\}\vert\partit)=\nicefrac{1}{3}I_{\{r,g\}}$ and therefore
\begin{equation*}
  \lpr(\lpr(\{g\}\vert\partit))
  =\frac{3}{4}\frac{\nicefrac{1}{3}+\nicefrac{1}{3}}{3}+\frac{1}{4}0
  =\frac{1}{6},
\end{equation*}
whereas $\lpr(\{g\})=\nicefrac{1}{4}$, and therefore $\lpr(\{g\})>\lpr(\lpr(\{g\}\vert\partit))$. $\blacklozenge$
\end{example}
\par
The difference $\upr(f\vert B)-\lpr(f\vert B)$ between infimum selling and supremum buying prices for gambles $f$ represents imprecision present in our subject's belief model. If we look at the inequalities in Proposition~\ref{prop:walley}.1, we are led to consider two extreme cases. One extreme maximises the `degrees of imprecision' $\upr(f\vert B)-\lpr(f\vert B)$ by letting $\lpr(f\vert B)=\inf_{\posty\in B}f(\posty)$ and $\upr(f\vert B)=\sup_{\posty\in   B}f(\posty)$. This leads to the so-called \emph{vacuous model}, corresponding to $\rdesirs=\set{f}{f\geq0}$, and intended to represent complete ignorance on the subject's part.
\par
The other extreme minimises the degrees of imprecision $\upr(f\vert B)-\lpr(f\vert B)$ by letting $\lpr(f\vert B)=\upr(f\vert B)$ everywhere. The common value $\pr(f\vert B)$ is then called the \emph{prevision}, or \emph{fair price}, for $f$ conditional on $B$. We call the corresponding functional $\pr(\cdot\vert B)$ a (conditional) \emph{linear prevision}. Linear previsions are the precise probability models considered by de Finetti \cite{finetti19745}. They of course have all properties of lower and upper previsions listed in Proposition~\ref{prop:walley}, with equality rather than inequality for statements 2 and 7. The restriction of a linear prevision to (indicators of) events is a finitely additive probability measure.

\section{Connecting the two approaches}\label{sec:connections}
In order to lay bare the connections between the game-theoretic and the behavioural approach, we enter Shafer and Vovk's world, and consider another player, called Forecaster, who, \emph{in situation $\init$}, has certain \emph{piece-wise} beliefs about what moves Reality will make.

\subsection{Forecaster's local beliefs}
More specifically, for each non-terminal situation $t\in\nonfins$, she has beliefs (in situation $\init$) about which move $\wmove$ Reality will choose from the set $\wmoves_t$ of moves available to him if he gets to $t$. We suppose she represents those beliefs in the form of a \emph{coherent}\footnote{Since we don't immediately envisage conditioning this local model on   subsets of $\wmoves_t$, we impose no extra conglomerability requirements here, only the   coherence conditions D1--D4.}  set $\rdesirs_t$ of really desirable gambles on $\wmoves_t$. These beliefs are conditional on the updating interpretation, in the sense that they represent Forecaster's beliefs in situation $\init$ about what Reality will do \emph{immediately after he   gets to situation $t$}. We call any specification of such coherent $\rdesirs_t$, $t\in\nonfins$, an \emph{immediate prediction model} for Forecaster. We want to stress here that $\rdesirs_t$ should \emph{not} be interpreted dynamically, i.e., as a set of gambles on $\wmoves_t$ that Forecaster accepts in situation $t$.
\par
We shall generally call an event tree, provided with local predictive belief models in each of the non-terminal situations $t$, an \emph{imprecise probability tree}. These local belief models may be coherent sets of really desirable gambles $\rdesirs_t$. But they can also be lower previsions $\lpr_t$ (perhaps derived from such sets $\rdesirs_t$). When all such local belief models are precise previsions, or equivalently (finitely additive) probability measures, we simply get a \emph{probability tree} in Shafer's \cite[Chapter~3]{shafer1996} sense.

\subsection{From local to global beliefs}
We can now ask ourselves what the behavioural implications of these conditional assessments $\rdesirs_t$ in the immediate prediction model are.  For instance, what do they tell us about whether or not Forecaster should accept certain gambles\footnote{In Shafer and Vovk's language, gambles are real variables.} on $\pspace$, the set of possible paths for Reality? In other words, how can these beliefs (in $\init$) about which next move Reality will make in each non-terminal situation $t$ be combined coherently into beliefs (in $\init$) about Reality's complete sequence of moves?
\par
In order to investigate this, we use Walley's very general and powerful method of \emph{natural extension}, which is just \emph{conservative coherent reasoning}. We shall construct, using the local pieces of information $\rdesirs_t$, a set of really desirable gambles on $\pspace$ for Forecaster in situation $\init$ that is (i) coherent, and (ii) as small as possible, meaning that no more gambles should be accepted than is actually required by coherence.

\subsubsection{Collecting the pieces}
\begin{figure}[h]
  \centering
  \begin{tikzpicture}[level distance=30mm]
    \tikzstyle{terminal}=[circle,draw=black!70,fill=black!70] 
    \tikzstyle{root}=[rectangle,draw=black!60]   
    \tikzstyle{nonterminal}=[circle,draw=black!60,fill=black!30,minimum size=.5mm] 
    \tikzstyle{level 1}=[sibling distance=15mm] 
    \tikzstyle{level 2}=[sibling distance=15mm] 
    \tikzstyle{level 3}=[sibling distance=15mm]
    \tikzstyle{curly}=[snake=snake,segment amplitude=.2mm,segment length=1mm, line after snake=1mm]
    \node[root] (root) {} [grow=right,inner sep=.7mm] 
    child {node[nonterminal,label=below:$t$] (t) {}
      child {node[nonterminal] (tw2) {}
        child {node[terminal] (fin1) {}}
        child {node[terminal] (fin2) {}}
        edge from parent {node[below] {$\wmove_2$}}}
      child {node[terminal] (tw1) {}
        edge from parent {node[above] {$\wmove_1$}}}}
    child {node[terminal] {} };
    \draw[darkgray,dashed] (t) -- +(-2,-1) node[below]    
    {$h_t\in\rdesirs_t\subseteq\gambles(\{\wmove_1,\wmove_2\})$}; 
    \draw[->,darkgray,curly] (fin1) -- +(1,0) node[right] {$h_t(\wmove_2)$};
    \draw[->,darkgray,curly] (fin2) -- +(1,0) node[right] {$h_t(\wmove_2)$};
    \draw[->,darkgray,curly] (tw1) -- +(1,0) node[right] {$h_t(\wmove_1)$};
    \draw[->,darkgray,curly] (tw2) -- +(0,-1) node[below] {$h_t(\wmove_2)$};
  \end{tikzpicture}
  \caption{In a non-terminal situation $t$, we consider a gamble $h_t$ on Reality's move space $\wmoves_t$ that Forecaster accepts, and turn it into a process, also denoted by $h_t$. The values $h_t(s)$ in situations $s\sfollows t$ are indicated by curly arrows.}
  \label{fig:local} 
\end{figure}
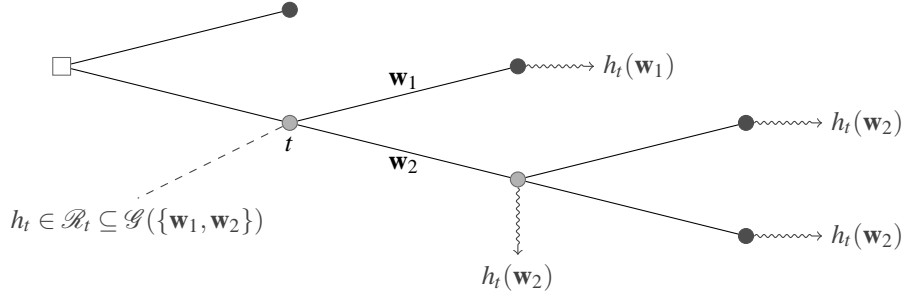

Consider any non-terminal situation $t\in\nonfins$ and any gamble $h_t$ in $\rdesirs_t$.  With $h_t$ we can associate a $t$-gamble,\footnote{Just as for variables, we can define a $t$-gamble as a partial gamble   whose domain includes $\upset{t}$.} also denoted by $h_t$, and defined by
\begin{equation*}
  h_t(\fin):= h_t(\fin(t))
\end{equation*}
for all $\fin\follows t$, where we denote by $\fin(t)$ the unique element of $\wmoves_t$ such that $t\fin(t)\precedes\fin$. The $t$-gamble $h_t$ is $U$-measurable for any cut $U$ of $t$ that is non-trivial, i.e., such that $U\neq\{t\}$. This implies that we can interpret $h_t$ as a map on $U$.  In fact, we shall even go further, and associate with the gamble $h_t$ on $\wmoves_t$ a $t$-process, also denoted by $h_t$, by letting $h_t(s):=h_t(\fin(t))$ for any $s\follows t$, where $\fin$ is any terminal situation that follows $s$; see also Fig.~\ref{fig:local}.
\par
\par
$I_{\upset{t}}h_t$ represents the gamble on $\fins$ that is called off unless Reality ends up in situation $t$, and which, when it isn't called off, depends only on Reality's move immediately after $t$, and gives the same value $h_t(\wmove)$ to all paths $\fin$ that go through $t\wmove$. The fact that Forecaster, in situation $\init$, accepts $h_t$ on $\wmoves_t$ conditional on Reality's getting to $t$, translates immediately to the fact that Forecaster accepts the contingent gamble $I_{\upset{t}}h_t$ on $\fins$, by Walley's Updating Principle.  We thus end up with a set 
\begin{equation*}
  \rdesirs
  :=\bigcup_{t\in\nonfins}\set{I_{\upset{t}}h_t}{h_t\in\rdesirs_t}
\end{equation*}
of gambles on $\fins$ that Forecaster accepts in situation $\init$. 
\par
The only thing left to do now, is to find the smallest coherent set $\natex_\rdesirs$ of really desirable gambles that includes $\rdesirs$ (if indeed there is any such coherent set).  Here we take coherence to refer to conditions D1--D4, together with D5', a variation on D5 which refers to conglomerability with respect to those partitions that we actually intend to condition on, as suggested in Section~\ref{sec:walley}.

\subsubsection{Cut conglomerability}
These partitions are what we call \emph{cut partitions}. Consider any cut $U$ of the initial situation $\init$. The set of events $\partit_U:=\set{\upset{u}}{u\in U}$ is a partition of $\fins$, called the \emph{$U$-partition}. D5' requires that our set of really desirable gambles should be \emph{cut conglomerable}, i.e., conglomerable with respect to every cut partition $\partit_U$.\footnote{Again, when all of Reality's move   spaces $\wmoves_t$ are finite, cut conglomerability (D5') is a consequence of D3, and therefore needs no extra   attention. But when some or all move spaces are infinite, then a cut $U$ may contain an infinite number of elements, and the corresponding cut partition $\partit_U$ will then be infinite too, making cut conglomerability a non-trivial additional requirement.}
\par
Why do we only require conglomerability for cut partitions?  Simply because we are interested in \emph{predictive inference}: we eventually will want to find out about the gambles on $\fins$ that Forecaster accepts in situation $\init$, conditional (contingent) on Reality getting to a situation $t$. This is related to finding lower previsions for Forecaster conditional on the corresponding events $\upset{t}$. A collection $\set{\upset{t}}{t\in T}$ of such events constitutes a partition of the sample space $\fins$ if and only if $T$ is a cut of $\init$.
\par
Because we require cut conglomerability, it follows in particular that $\natex_\rdesirs$ will contain the sums of gambles $g:=\sum_{u\in U}I_{\upset{u}}h_u$ for all \emph{non-terminal} cuts $U$ of $\square$ and all choices of $h_u\in\rdesirs_u$, $u\in U$.  This is because $I_{\upset{u}}g=I_{\upset{u}}h_u\in\rdesirs$ for all $u\in U$. Because moreover $\natex_\rdesirs$ should be a convex cone [by D3 and D4], any sum of such sums $\sum_{u\in U}I_{\upset{u}}h_u$ over a finite number of non-terminal cuts $U$ should also belong to $\natex_\rdesirs$. But, since in the case of bounded protocols we are discussing here, Reality can only make a bounded and finite number of moves, $\nonfins$ is a finite union of such non-terminal cuts, and therefore the sums $\sum_{u\in\nonfins}I_{\upset{u}}h_u$ should belong to $\natex_\rdesirs$ for all choices $h_u\in\rdesirs_u$, $u\in\nonfins$.

\subsubsection{Selections and gamble processes}
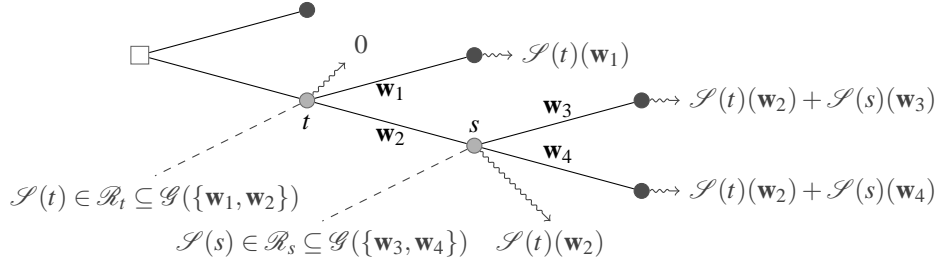
\begin{figure}[h]
  \centering
  \begin{tikzpicture}[level distance=22mm]
    \tikzstyle{terminal}=[circle,draw=black!70,fill=black!70] 
    \tikzstyle{root}=[rectangle,draw=black!60]   
    \tikzstyle{nonterminal}=[circle,draw=black!60,fill=black!30,minimum size=.5mm] 
    \tikzstyle{level 1}=[sibling distance=12mm] 
    \tikzstyle{level 2}=[sibling distance=12mm] 
    \tikzstyle{level 3}=[sibling distance=12mm]
    \tikzstyle{curly}=[snake=snake,segment amplitude=.2mm,segment length=1mm, line after snake=1mm]
    \node[root] (root) {} [grow=right,inner sep=.7mm] 
    child {node[nonterminal,label=below:$t$] (t) {}
      child {node[nonterminal,label=above:$s$] (s) {}
        child {node[terminal] (fin1) {}
          edge from parent {node[above] {$\wmove_4$}}}
        child {node[terminal] (fin2) {}
          edge from parent {node[above] {$\wmove_3$}}}
        edge from parent {node[below] {$\wmove_2$}}}
      child {node[terminal] (tw1) {}
        edge from parent {node[below] {$\wmove_1$}}}}
    child {node[terminal] {} };
    \draw[darkgray,dashed] (t) -- +(-2,-1) node[below]    
    {$\selec(t)\in\rdesirs_t\subseteq\gambles(\{\wmove_1,\wmove_2\})$}; 
    \draw[darkgray,dashed] (s) -- +(-2,-1) node[below]    
    {$\selec(s)\in\rdesirs_s\subseteq\gambles(\{\wmove_3,\wmove_4\})$}; 
    \draw[->,darkgray,curly] (fin1) -- +(.5,0) node[right] {$\selec(t)(\wmove_2)+\selec(s)(\wmove_4)$};
    \draw[->,darkgray,curly] (fin2) -- +(.5,0) node[right] {$\selec(t)(\wmove_2)+\selec(s)(\wmove_3)$};
    \draw[->,darkgray,curly] (tw1) -- +(.5,0) node[right] {$\selec(t)(\wmove_1)$};
    \draw[->,darkgray,curly] (s) -- +(1,-1) node[below] {$\selec(t)(\wmove_2)$};
    \draw[->,darkgray,curly] (t) -- +(.5,.5) node[above right] {$0$};
  \end{tikzpicture}
  \caption{The $t$-selection $\selec$ in this event tree is a process defined in the two non-terminal situations $t$ and $s$; it selects, in each of these situations, a really desirable gamble for Forecaster. The values of the corresponding gamble process $\gamble^\selec$ are indicated by curly arrows.}
  \label{fig:gamble-processes} 
\end{figure}

Consider any  non-terminal situation $t$, and call \emph{$t$-selection} any partial process $\selec$ defined on the non-terminal $s\follows t$ such that $\selec(s)\in\rdesirs_s$. With a $t$-selection $\selec$, we associate a $t$-process $\gamble^\selec$, called a \emph{gamble process}, where
\begin{equation}\label{eq:gamble-process}
  \gamble^\selec(s)
  =\sum_{t\precedes u\sprecedes s}\selec(u)(s)
\end{equation}
in all situations $s\follows t$; see also Fig.~\ref{fig:gamble-processes}. Alternatively, $\gamble^\selec$ is given by the recursion relation
\begin{equation*}
  \gamble^\selec(s\wmove) =\gamble^\selec(s)+\selec(s)(\wmove),\quad\wmove\in\wmoves_s
\end{equation*}
for all non-terminal $s\follows t$, with initial value $\gamble^\selec(t)=0$. In particular, this leads to the $t$-gamble $\gamble^\selec_\fins$ defined on all terminal situations $\fin$ that follow $t$, by letting
\begin{equation}\label{eq:gamble-variable}
  \gamble^\selec_\fins
  =\sum_{t\precedes u,u\in\nonfins}I_{\upset{u}}\selec(u).
\end{equation}
Then we have just argued that the gambles $\gamble^\selec_\fins$ should belong to $\natex_\rdesirs$ for all non-terminal situations $t$ and all $t$-selections $\selec$. As before for strategy and capital processes, we call a $\init$-selection $\selec$ simply a \emph{selection}, and a $\init$-gamble process simply a \emph{gamble process}. 
\par

\subsubsection{The Marginal Extension Theorem}
It is now but a technical step to prove Theorem~\ref{theo:natex} below. It is a significant generalisation, in terms of sets of really desirable gambles rather than coherent lower previsions,\footnote{The difference in language may obscure that this is indeed a generalisation.   But see Theorem~\ref{theo:concatenation} for expressions in terms of predictive lower previsions   that should make the connection much clearer.} of the \emph{Marginal Extension Theorem} first proved by Walley \cite[Theorem~6.7.2]{walley1991}, and subsequently extended by De Cooman and Miranda \cite{miranda2006b}.

\begin{theorem}[Marginal Extension Theorem]\label{theo:natex}
  There is a smallest set of gambles that satisfies D1--D4 and D5' and includes $\rdesirs$. This   natural extension of\/ $\rdesirs$ is given by
\begin{equation*}
  \natex_\rdesirs
  :=\set{g}
  {\text{$g\geq\gamble^\selec_\fins$ for some selection $\selec$}}.
\end{equation*}
Moreover, for any non-terminal situation $t$ and any $t$-gamble $g$, it holds that\/ $I_{\upset{t}}g\in\natex_\rdesirs$ if and only if there is some $t$-selection $\selec_t$ such that $g\geq\gamble^{\selec_t}_\fins$, where as before, $g\geq\gamble^{\selec_t}_\fins$ is taken to mean that $g(\fin)\geq\gamble^{\selec_t}_\fins(\fin)$ for all terminal situations $\fin$ that follow $t$.
\end{theorem}

\subsection{Predictive lower and upper previsions}\label{sec:predictive-lower-upper}
We now use the coherent set of really desirable gambles $\natex_\rdesirs$ to define special lower previsions $\lpr(\cdot\vert t):=\lpr(\cdot\vert\upset{t})$ for Forecaster in situation $\init$, conditional on an event $\upset{t}$, i.e., on Reality getting to situation $t$, as explained in Section~\ref{sec:walley}.\footnote{We stress   again that these are conditional lower previsions on the contingent/updating interpretation.} We shall call such conditional lower previsions \emph{predictive} lower previsions. We then get, using Eq.~\eqref{eq:lowdef} and Theorem~\ref{theo:natex}, that for any non-terminal situation $t$,
\begin{align}
  \lpr(f\vert t) &:=\sup\set{\alpha}
  {I_{\upset{t}}(f-\alpha)\in\natex_\rdesirs}\label{eq:lprdef}\\
  &=\sup\set{\alpha} {\text{$f-\alpha\geq\gamble^\selec_\fins$ for some $t$-selection $\selec$}}.\label{eq:lpr}
\end{align}
We also use the notation $\lpr(f):=\lpr(f\vert\init)=\sup\set{\alpha}{f-\alpha\in\natex_\rdesirs}$. It should be stressed that Eq.~\eqref{eq:lprdef} is also valid in terminal situations $t$, whereas Eq.~\eqref{eq:lpr} clearly isn't.
\par
Besides the properties in Proposition~\ref{prop:walley}, which hold in general for conditional lower and upper previsions, the predictive lower (and upper) previsions we consider here also satisfy a number of additional properties, listed in Propositions~\ref{prop:prevision-properties-general} and~\ref{prop:separate-coherence}.

\begin{proposition}[Additional properties of predictive lower and upper
  previsions]\label{prop:prevision-properties-general}
  Let\/ $t$ be any situation, and let\/ $f$, $f_1$ and $f_2$ be gambles on $\fins$.
  \begin{enumerate}[1.]
  \item If\/ $t$ is a terminal situation $\fin$, then $\lpr(f\vert\fin)=\upr(f\vert\fin)=f(\fin)$;
  \item $\lpr(f\vert t)=\lpr(fI_{\upset{t}}\vert t)$ and $\upr(f\vert t)=\upr(fI_{\upset{t}}\vert t)$;
  \item $f_1\leq f_2$ (on $\upset{t}$) implies that $\lpr(f_1\vert t)\leq\lpr(f_2\vert t)$ [monotonicity].
  \end{enumerate}
\end{proposition}
\noindent
Before we go on, there is an important point that must be stressed and clarified. It is an immediate consequence of Proposition~\ref{prop:prevision-properties-general}.2 that when $f$ and $g$ are any two gambles that coincide on $\upset{t}$, then $\lpr(f\vert t)=\lpr(g\vert t)$. This means that $\lpr(f\vert t)$ is completely determined by the values that $f$ assumes on $\upset{t}$, and it allows us to define $\lpr(\cdot\vert t)$ on gambles that are only necessarily defined on $\upset{t}$, i.e., on $t$-gambles. We shall do so freely in what follows.
\par
For any cut $U$ of a situation $t$, we may define the $t$-gamble $\lpr(f\vert U)$ as the gamble that assumes the value $\lpr(f\vert u)$ in any $\fin\follows u$, where $u\in U$. This $t$-gamble is $U$-measurable by construction, and it can be considered as a gamble on $U$.
\begin{proposition}[Separate coherence]\label{prop:separate-coherence}
  Let\/ $t$ be any situation, let\/ $U$ be any cut of\/ $t$, and let\/ $f$ and $g$ be $t$-gambles, where $g$ is   $U$-measurable.
  \begin{enumerate}[1.]
  \item $\lpr(\upset{t}\vert t)=1$;
  \item $\lpr(g\vert U)=g_U$;
  \item $\lpr(f+g\vert U)=g_U+\lpr(f\vert U)$;
  \item if $g$ is moreover non-negative, then $\lpr(gf\vert U)=g_U\lpr(f\vert U)$.
  \end{enumerate}
\end{proposition}

\subsection{Correspondence between immediate prediction models and coherent probability protocols}
There appears to be a close correspondence between the expressions [such as~\eqref{eq:lprice}] for lower prices $\lprice_t(f)$ associated with coherent probability protocols and those [such as~\eqref{eq:lpr}] for the predictive lower previsions $\lpr(f\vert t)$ based on an immediate prediction model.  Say that a given coherent probability protocol and given immediate prediction model \emph{match} whenever they lead to identical corresponding lower prices $\lprice_t$ and predictive lower previsions $\lpr(\cdot\vert t)$ for all \emph{non-terminal} $t\in\nonfins$.
\par
The following theorem marks the culmination of our search for the correspondence between Walley's, and Shafer and Vovk's approaches to probability theory.

\begin{theorem}[Matching Theorem]\label{theo:matching}
  For every coherent probability protocol there is an immediate prediction model such that the two   match, and conversely, for every immediate prediction model there is a coherent probability protocol   such that the two match.
\end{theorem}
\noindent
The ideas underlying the proof of this theorem should be clear. If we have a coherent probability protocol with move spaces $\smoves_t$ and gain functions $\gain_t$ for Sceptic, define the immediate prediction model for Forecaster to be (essentially) $\rdesirs_t:=\set{-\gain(\smove,\cdot)}{\smove\in\smoves_t}$.  If, conversely, we have an immediate prediction model for Forecaster consisting of the sets $\rdesirs_t$, define the move spaces for Sceptic by $\smoves_t:=\rdesirs_t$, and his gain functions by $\gain_t(h,\cdot):=-h$ for all $h$ in $\rdesirs_t$. We discuss the interpretation of this correspondence in more detail in Section~\ref{sec:interpretation}. 

\subsection{Calculating predictive lower prevision using backwards recursion}\label{sec:concatenation}
The Marginal Extension Theorem allows us to calculate the most conservative global belief model $\natex_\rdesirs$ that corresponds to the local immediate prediction models $\rdesirs_t$. Here beliefs are expressed in terms of sets of really desirable gambles. Can we derive a result that allows us to do something similar for the corresponding lower previsions?
\par
To see what this question entails, first consider a local model $\rdesirs_s$: a set of really desirable gambles on $\wmoves_s$, where $s\in\nonfins$. Using Eq.~\eqref{eq:lowdef}, we can associate with $\rdesirs_s$ a lower prevision $\lpr_s$ on $\gambles(\wmoves_s)$. Each gamble $g_s$ on $\wmoves_s$ can be seen as an uncertain reward, whose outcome $g_s(\wmove)$ depends on the (unknown) move $\wmove\in\wmoves_s$ that Reality will make if it gets to situation $s$. And Forecaster's \emph{local} (predictive) lower prevision 
\begin{equation}\label{eq:local-prevision}
  \lpr_s(g_s):=\sup\set{\alpha}{g_s-\alpha\in\rdesirs_s}
\end{equation}
for $g_s$ is her supremum acceptable price (in $\init$) for buying $g_s$ when Reality gets to $s$. 
\par
But as we have seen in Section~\ref{sec:predictive-lower-upper}, we can also, in each situation $t$, derive \emph{global} predictive lower previsions $\lpr(\cdot\vert t)$ for Forecaster from the global model $\natex_\rdesirs$, using Eq.~\eqref{eq:lprdef}. For each $t$-gamble $f$, $\lpr(f\vert t)$ is Forecaster inferred supremum acceptable price (in $\init$) for buying $f$, contingent on Reality getting to $t$. 
\par
Is there a way to construct the global predictive lower previsions $\lpr(\cdot\vert t)$ directly from the local predictive lower previsions $\lpr_s$? We can infer that there is from the following theorem, together with Propositions~\ref{prop:cut-reduction} and~\ref{prop:local-models} below.

\begin{theorem}[Concatenation Formula]\label{theo:concatenation}
  Consider any two cuts $U$ and $V$ of a situation $t$ such that\/ $U$ precedes $V$. For all\/   $t$-gambles $f$ on $\fins$,\footnote{Here too, it is implicitly assumed that all expressions     are well-defined, e.g., that in the second statement, $\lpr(f\vert v)$ is a real number for     all $v\in V$, making sure that $\lpr(f\vert V)$ is indeed a gamble.}
  \begin{enumerate}[1.]
  \item $\lpr(f\vert t)=\lpr(\lpr(f\vert U)\vert t)$;
  \item $\lpr(f\vert U)=\lpr(\lpr(f\vert V)\vert U)$.
  \end{enumerate}
\end{theorem}
\noindent
To make clear what the following Proposition~\ref{prop:cut-reduction} implies, consider any $t$-selection $\selec$, and define the \emph{$U$-called off $t$-selection} $\selec^U$ as the selection that mimics $\selec$ until we get to $U$, where we begin to select the zero gambles: for any non-terminal situation $s\follows t$, let $\selec^U(s):=\selec(s)$ if $s$ strictly precedes (some element of) $U$, and let $\selec^U(s):=0\in\rdesirs_s$ otherwise. If we stop the gamble process $\gamble^\selec$ at the cut $U$, we readily infer from Eq.~\eqref{eq:gamble-process} that for the $U$-stopped process $U(\gamble^\selec)$
\begin{equation}\label{eq:stopped-too}
  U(\gamble^\selec)=\gamble^{\selec^U}
  \text{ and therefore,  also using Eq.~\eqref{eq:stopped}, }
  \gamble^\selec_U=\gamble^{\selec^U}_\fins.
\end{equation}
We see that stopped gamble processes are gamble processes themselves, that correspond to selections being `called off' as soon as Reality reaches a cut. This also means that we can actually restrict ourselves to selections $\selec$ that are $U$-called off in Proposition~\ref{prop:cut-reduction}.

\begin{proposition}\label{prop:cut-reduction}
  Let\/ $t$ be a non-terminal situation, and let\/ $U$ be a cut of\/ $t$. Then for any   $U$-measurable $t$-gamble $f$, $I_{\upset{t}}f\in\natex_\rdesirs$ if and only is there is some   $t$-selection $\selec$ such that\/ $I_{\upset{t}}f\geq\gamble^{\selec^U}_\fins$, or   equivalently, $f_U\geq\gamble^\selec_U$. Consequently,
  \begin{align*}
    \lpr(f\vert t) &=\sup\set{\alpha}{\text{$f-\alpha\geq\gamble^{\selec^U}_\fins$ for some
        $t$-selection $\selec$}}\\
    &=\sup\set{\alpha}{\text{$f_U-\alpha\geq\gamble^{\selec}_U$ for some $t$-selection $\selec$}}.
  \end{align*}
\end{proposition}
\noindent
If a $t$-gamble $h$ is measurable with respect to the children cut $C(t)$ of a non-terminal situation $t$, then we can interpret it as gamble on $\wmoves_t$.  For such gambles, the following immediate corollary of Proposition~\ref{prop:cut-reduction} tells us that the predictive lower previsions $\lpr(h\vert t)$ are completely determined by the local modal $\rdesirs_t$.

\begin{proposition}\label{prop:local-models}
  Let\/ $t$ be a non-terminal situation, and consider a $C(t)$-measurable gamble $h$. Then $\lpr(h\vert t)=\lpr_t(h)$.
\end{proposition}

These results tells us that all predictive lower (and upper) previsions can be calculated using backwards recursion, by starting with the trivial predictive previsions $\upr(f\vert\fins)=\lpr(f\vert\fins)=f$ for the terminal cut $\fins$, and using only the local models $\lpr_t$. This is illustrated in the following simple example. We shall come back to this idea in Section~\ref{sec:backwards-recursion}.

\begin{example}\label{ex:many-coins}
  Suppose we have $n>0$ coins. We begin by flipping the first coin: if we get tails, we stop, and otherwise we flip   the second coin. Again, we stop if we get tails, and otherwise we flip the third coin, \dots\ In other words, we   continue flipping new coins until we get one tails, or until all $n$ coins have been flipped. This leads to the   event tree depicted in Fig.~\ref{fig:many-coins}.
\begin{figure}[h]
  \centering
  \begin{tikzpicture}[level distance=20mm,sibling distance=10mm]
    \tikzstyle{terminal}=[circle,draw=black!70,fill=black!70] 
    \tikzstyle{root}=[rectangle,draw=black!60]   
    \tikzstyle{nonterminal}=[circle,draw=black!60,fill=black!30,minimum size=.5mm] 
    \tikzstyle{cut}=[semithick,dashed]
    \node[root,label=below:$\heads_0$] (root) {} [grow=right,inner sep=.7mm] 
    child {node[nonterminal,label=below left:$\heads_1$] (h1) {}
      child {node[nonterminal,label=below left:$\heads_2$] (h2) {}
        child {node[nonterminal,label=below left:$\heads_3$] (h3) {}
          child {node[nonterminal,label=below left:$\heads_{n-1}$] (hnm1) {}
            child {node[terminal,label=right:$\heads_n$] (hn) {} edge from parent[solid]}
            child {node[terminal,label=right:$\tails_n$] (tn) {} edge from parent[solid]}
            edge from parent[dotted]}
          child {node[terminal,label=right:$\tails_4$] (t4) {}}}
        child {node[terminal,label=right:$\tails_3$] (t3) {}}}
      child {node[terminal,label=right:$\tails_2$] (t2) {}}}
    child {node[terminal,label=right:$\tails_1$] (t1) {}};
    \draw[cut] (t1) -- +(0,.5);
    \draw[cut] (t1) -- (h1) -- +(0,-1) node[right] {$U_1$};
    \draw[cut] (t2) -- +(0,.5);
    \draw[cut] (t2) -- (h2) -- +(0,-1) node[right] {$U_2$};
    \draw[cut] (t3) -- +(0,.5);
    \draw[cut] (t3) -- (h3) -- +(0,-1) node[right] {$U_3$};
    \draw[cut] (hnm1) -- +(0,.5);
    \draw[cut] (hnm1) -- +(0,-1) node[right] {$U_{n-1}$};
    \draw[cut] (tn) -- +(0,.5);
    \draw[cut] (tn) -- (hn) -- +(0,-1) node[right] {$U_n$};
  \end{tikzpicture}
  \caption{The event tree for the uncertain process involving $n$ successive coin flips described in Example~\ref{ex:many-coins}.}  
  \label{fig:many-coins}
\end{figure}
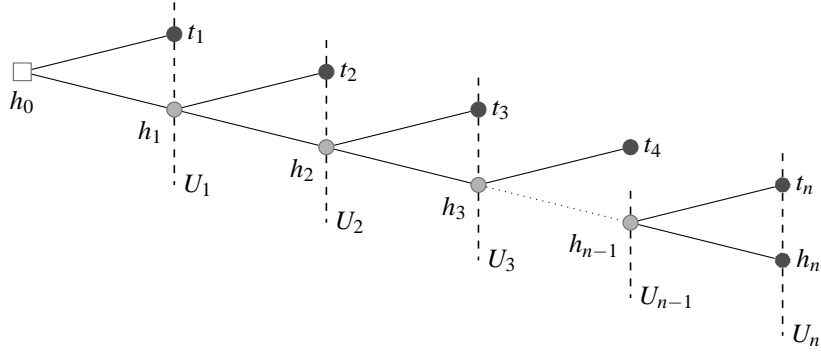
Its sample space is $\fins=\{\tails_1,\tails_2,\dots,\tails_n,\heads_n\}$. We will also consider the cuts $U_1=\{\tails_1,\heads_1\}$ of $\init$, $U_2=\{\tails_2,\heads_2\}$ of $\heads_1$, $U_3=\{\tails_3,\heads_3\}$ of $\heads_2$, \dots, and $U_n=\{\tails_n,\heads_n\}$ of $\heads_{n-1}$. It will be convenient to also introduce the notation $\heads_0$ for the initial situation $\init$.
\par
For each of the non-terminal situations $\heads_k$, $k=0,1,\dots,n-1$, Forecaster has beliefs (in $\init$) about what move Reality will make in that situation, i.e., about the outcome of the $k+1$-th coin flip. These beliefs are expressed in terms of a set of really desirable gambles $\rdesirs_{\heads_k}$ on Reality's move space $\wmoves_{\heads_k}$ in $\heads_k$. Each such move space $\wmoves_{\heads_k}$ can clearly be identified with the children cut $U_{k+1}$ of $\heads_k$.
\par
For the purpose of this example, it will be enough to consider the local predictive lower previsions $\lpr_{\heads_k}$ on $\gambles(U_{k+1})$, associated with $\rdesirs_{\heads_k}$ through Eq.~\eqref{eq:local-prevision}.  Forecaster assumes all coins to be approximately fair, in the sense that she assesses that the probability of heads for each flip lies between $\frac{1}{2}-\delta$ and $\frac{1}{2}+\delta$, for some $0<\delta<\frac{1}{2}$. This assessment leads to the following local predictive lower previsions:\footnote{These so-called linear-vacuous mixtures, or contamination   models, are the natural extensions of the probability assessments   $\lpr_{\heads_k}(\{\heads_{k+1}\})=\frac{1}{2}-\delta$ and   $\upr_{\heads_k}(\{\heads_{k+1}\})=\frac{1}{2}+\delta$; see \cite[Chapters~3--4]{walley1991} for more details.}
\begin{equation}\label{eq:stopped-flips}
  \lpr_{\heads_k}(g)
  =(1-2\delta)\left[\frac{1}{2}g(\heads_{k+1})+\frac{1}{2}g(\tails_{k+1})\right]
  +2\delta\min\{g(\heads_{k+1}),g(\tails_{k+1})\},
\end{equation}
where $g$ is any gamble on $U_{k+1}$.
\par
Let us see how we can for instance calculate, from the local predictive models $\lpr_{\heads_k}$, the predictive lower probabilities $\lpr(\{\heads_n\}\vert s)$ for a gamble $f$ on $\fins$ and any situation $s$ in the tree.  First of all, for the terminal situations it is clear from Proposition~\ref{prop:prevision-properties-general}.1 that
\begin{equation}\label{eq:stopped-flips-1}
  \lpr(\{\heads_n\}\vert\tails_k)=0\text{ and }\lpr(\{\heads_n\}\vert\heads_n)=1.
\end{equation}
We now turn to the calculation of $\lpr(\{\heads_n\}\vert\heads_{n-1})$. It follows at once from Proposition~\ref{prop:local-models} that $\lpr(\{\heads_n\}\vert\heads_{n-1})=\lpr_{\heads_{n-1}}(\{\heads_n\})$, and therefore, substituting $g=I_{\{\heads_n\}}$ in Eq.~\eqref{eq:stopped-flips} for $k=n-1$,
\begin{equation}\label{eq:stopped-flips-2}
  \lpr(\{\heads_n\}\vert\heads_{n-1})=\frac{1}{2}-\delta.
\end{equation}
To calculate $\lpr(\{\heads_n\}\vert\heads_{n-2})$, consider that, since $h_{n-1}\precedes U_{n-1}$,
\begin{equation*}
  \lpr(\{\heads_n\}\vert\heads_{n-2})
  =\lpr(\lpr(\{\heads_n\}\vert U_{n-1})\vert\heads_{n-2})
  =\lpr_{\heads_{n-2}}(\lpr(\{\heads_n\}\vert U_{n-1}))
\end{equation*}
where the first equality follows from Theorem~\ref{theo:concatenation}, and the second from Proposition~\ref{prop:local-models}, taking into account that $g_{n-1}:=\lpr(\{\heads_n\}\vert U_{n-1})$ is a gamble on the children cut $U_{n-1}$ of $\heads_{n-2}$. It follows from Eq.~\eqref{eq:stopped-flips-1} that $g_{n-1}(\tails_{n-1})=\lpr(\{\heads_n\}\vert\tails_{n-1})=0$ and from Eq.~\eqref{eq:stopped-flips-2} that $g_{n-1}(\heads_{n-1})=\lpr(\{\heads_n\}\vert\heads_{n-1})=\frac{1}{2}-\delta$. Substituting $g=g_{n-1}$ in Eq.~\eqref{eq:stopped-flips} for $k=n-2$, we then find that
\begin{equation}\label{eq:stopped-flips-3}
  \lpr(\{\heads_n\}\vert\heads_{n-2})=(\frac{1}{2}-\delta)^2.
\end{equation}
Repeating this course of reasoning, we find that more generally
\begin{equation}\label{eq:stopped-flips-4}
  \lpr(\{\heads_n\}\vert\heads_{k})=(\frac{1}{2}-\delta)^{n-k},\quad k=0,\dots n-1.
\end{equation}
This illustrates how we can use a backwards recursion procedure to calculate global from local predictive lower previsions.\footnote{It also indicates why we need to work in the more general language of lower previsions and gambles, rather than the perhaps more familiar one of lower probabilities and events: even if we only want to calculate a global predictive lower probability, already after one recursion step we need to start working with lower previsions of gambles. More discussion on the prevision/gamble versus probability/event issue can be found in \cite[Chapter~4]{walley1991}.} 
\end{example}

\section{Interpretation of the Matching Theorem}\label{sec:interpretation}
In Shafer and Vovk's approach, there sometimes also appears, besides Reality and Sceptic, a third player, called \emph{Forecaster}. Her r\^ole consists in determining what Sceptic's move space $\smoves_t$ and gain function $\gain_t$ are, in each non-terminal situation $t$. Shafer and Vovk leave largely unspecified just how Forecaster should do that, which makes their approach quite general and abstract.
\par
But the Matching Theorem now tells us that we can connect their approach with Walley's, and therefore inject a notion of belief modelling into their game-theoretic framework. We can do that by being more specific about how Forecaster should determine Sceptic's move spaces $\smoves_t$ and gain functions $\gain_t$: they should be determined by Forecaster's beliefs (in $\init$) about what Reality will do immediately after getting to non-terminal situations $t$.\footnote{The germ for this idea, in the case that Forecaster's beliefs can be expressed using precise   probability models on the $\gambles(\wmoves_t)$, is already present in Shafer's work, see for instance   \cite[Chapter~8]{shafer2001} and \cite[Appendix~1]{shafer1985}. We extend this idea here to Walley's imprecise   probability models.} Let us explain this more carefully.
\par
Suppose that Forecaster has certain beliefs, \emph{in situation $\init$}, about what move Reality will make next in each non-terminal situation $t$, and suppose she models those beliefs by specifying a coherent set $\rdesirs_t$ of really desirable gambles on $\wmoves_t$. This brings us to the situation described in the previous section. 
\par
When Forecaster specifies such a set, she is making certain behavioural commitments: she is committing herself to accepting, in situation $\init$, any gamble in $\rdesirs_t$, contingent on Reality getting to situation $t$, and to accepting any combination of such gambles according to the combination axioms D3, D4 and D5'. This implies that we can derive predictive lower previsions $\lpr(\cdot\vert t)$, with the following interpretation: in situation $\init$, $\lpr(f\vert t)$ is the supremum price Forecaster can be made to buy the $t$-gamble $f$ for, conditional on Reality's getting to $t$, and on the basis of the commitments she has made in the initial situation $\init$.
\par
What Sceptic can now do, is take Forecaster up on her commitments. This means that in situation $\init$, he can use a selection $\selec$, which for each non-terminal situation $t$, selects a gamble (or equivalently, any non-negative linear combination of gambles) $\selec(t)=h_t$ in $\rdesirs_t$ and offer the corresponding gamble $\gamble^\selec_\fins$ on $\fins$ to Forecaster, who is bound to accept it. If Reality's next move in situation $t$ is $\wmove\in\wmoves_t$, this changes Sceptic's capital by (the positive or negative amount) $-h_t(\wmove)$. In other words, his move space $\smoves_t$ can then be identified with the convex set of gambles $\rdesirs_t$ and his gain function $\gain_t$ is then given by $\gain_t(h_t,\cdot)=-h_t$. But then the \emph{selection} $\selec$ can be identified with a \emph{strategy} $\strat$ for Sceptic, and $\cptl^\strat_\fins=-\gamble^\selec_\fins$ (this is the essence of the proof of Theorem~\ref{theo:matching}), which tells us that we are led to a coherent probability protocol, and that the corresponding lower prices $\lprice_t$ for Sceptic coincide with Forecaster's predictive lower previsions $\lpr(\cdot\vert t)$.
\par
In a very nice paper \cite{shafer2003}, Shafer, Gillett and Scherl discuss ways of introducing and interpreting lower previsions in a game-theoretic framework, not in terms of prices that a subject is willing to pay for a gamble, but in terms of whether a subject believes she can make a lot of money (utility) at those prices. They consider such conditional lower previsions both on a contingent and on a dynamic interpretation, and argue that there is equality between them in certain cases. Here, we have decided to stick to the more usual interpretation of lower and upper previsions, and concentrated on the contingent/updating interpretation. We see that on our approach, the game-theoretic framework is useful too.
\par
This is of particular relevance to the laws of large numbers that Shafer and Vovk derive in their game-theoretic framework, because such laws can now be given a behavioural interpretation in terms of Forecaster's predictive lower and upper previsions. To give an example, we now turn to deriving a very general weak law of large numbers.

\section{A more general weak law of large numbers}
\label{sec:weak-law}
Consider a non-terminal situation $t$ and a cut $U$ of $t$.  Define the $t$-variable $n_U$ such that $n_U(\fin)$ is the distance $d(t,u)$, measured in moves along the tree, from $t$ to the unique situation $u$ in $U$ that $\fin$ goes through. $n_U$ is clearly $U$-measurable, and $n_U(u)$ is simply the distance $d(t,u)$ from $t$ to $u$. We assume that $n_U(u)>0$ for all $u\in U$, or in other words that $U\neq\{t\}$. Of course, in the bounded protocols we are considering here, $n_U$ is bounded, and we denote its minimum by $N_U$.
\par
Now consider for each $s$ between $t$ and $U$ a \emph{bounded} gamble $h_s$ and a real number $m_s$ such that $h_s-m_s\in\rdesirs_s$, meaning that Forecaster in situation $\init$ accepts to buy $h_s$ for $m_s$, contingent on Reality getting to situation $s$. Let $B>0$ be any common upper bound for $\sup h_s-\inf h_s$, for all $t\precedes s\sprecedes U$. It follows from the coherence of $\rdesirs_s$ [D1] that $m_s\leq\sup h_s$. To make things interesting, we shall also assume that $\inf h_s\leq m_s$, because otherwise $h_s-m_s\geq0$ and accepting this gamble represents no real commitment on Forecaster's part. As a result, we see that $\abs{h_s-m_s}\leq\sup h_s-\inf h_s\leq B$.
\par
We are interested in the following $t$-gamble $G_U$, given by
\begin{equation*}
  G_U=\frac{1}{n_U}\sum_{t\precedes s\sprecedes U}I_{\upset{s}}[h_s-m_s],
\end{equation*}
which provides a measure for how much, on average, the gambles $h_s$ yield an outcome above Forecaster's accepted buying prices $m_s$, along segments of the tree starting in $t$ and ending right before $U$. In other words, $G_U$ measures the average gain for Forecaster along segments from $t$ to $U$, associated with commitments she has made and is taken up on, because Reality has to move along these segments. This gamble $G_U$ is $U$-measurable too. We may therefore interpret $G_U$ as a gamble on $U$.  Also, for any $h_s$ and any $u\in U$, we know that because $s\sprecedes u$, $h_s$ has the same value $h_s(u):=h_s(\fin(s))$ in all $\fin$ that go through $u$. This allows us to write
\begin{equation*}
  G_U(u)=\frac{1}{n_U(u)}\sum_{t\precedes s\sprecedes u}[h_s(u)-m_s].
\end{equation*}
We would like to study Forecaster's beliefs (in the initial situation $\init$ and contingent on Reality getting to $t$) in the occurrence of the event
\begin{equation*}
  \propset{G_U\geq-\epsilon}:=\set{\fin\in\upset{t}}{G_U(\fin)\geq-\epsilon},
\end{equation*}
where $\epsilon>0$. In other words, we want to know $\lpr(\propset{G_U\geq-\epsilon}\vert t)$, which is Forecaster's supremum rate for betting on the event that his average gain from $t$ to $U$ will be at least $-\epsilon$, contingent on Reality's getting to $t$.

\begin{theorem}[Weak Law of Large Numbers]\label{theo:largenum}
  For all $\epsilon>0$,
  \begin{equation*}
    \lpr(\propset{G_U\geq-\epsilon}\vert t)
    \geq1-\exp\left(-\frac{N_U\epsilon^2}{4B^2}\right).
  \end{equation*}
\end{theorem}
\noindent
We see that as $N_U$ increases this lower bound increases to one, so the theorem can be very loosely formulated as follows: \emph{As the horizon recedes, Forecaster, if she is coherent, should believe increasingly more strongly that her average gain   along any path from the present to the horizon won't be negative}. 
\par
This is a very general version of the weak law of large numbers. It can be seen as a generalisation of Hoeffding's inequality for martingale differences \cite{hoeffding1963} (see also \cite[Chapter~4]{wasserman2004} and \cite[Appendix~A.7]{vovk2005}) to coherent lower previsions on event trees.

\section{Scoring a predictive model}\label{sec:scoring}
We now look at an interesting consequence of Theorem~\ref{theo:largenum}: we shall see that it can be used to score a predictive model in a manner that satisfies Dawid's \emph{Prequential Principle} \cite{dawid1984,dawid1999}. We consider the special case of  Theorem~\ref{theo:largenum} where $t=\init$.
\par
\begin{figure}[h]
  \centering
  \begin{tikzpicture}[samples=200,scale=5,domain=0:1]
    \draw[ultra thin,color=gray,step=.1] (-0.01,-0.01) grid (1.01,1.01); 
    \draw (1.0,0.0) node[below] {1}; 
    \draw (0.0,1.0) node[left] {1}; 
    \draw (0.0,0.0) node[below right] {0}; 
    \draw (0.0,0.0) node[above left] {0}; 
    \draw[->] (-0.05,0) -- (1.05,0) node[right] {$-x$}; 
    \draw[->] (0,-0.05) -- (0,1.05) node[above] {$1-S_N(x)$}; 
    \draw[smooth,color=blue] plot[id=5] function{1-exp(-x**2*(5/4))} node[right]{$N_U=5$}; 
    \draw[smooth,color=red] plot[id=10] function{1-exp(-x**2*(10/4))} node[right]{$N_U=10$}; 
    \draw[smooth,color=orange] plot[id=100] function{1-exp(-x**2*(100/4))} node[right]{$N_U=100$}; 
    \draw[smooth,color=brown] plot[id=500] function{1-exp(-x**2*(500/4))} node[above left]{$N_U=500$};
  \end{tikzpicture}
  \caption{What Forecaster can be made to pay, $1-S_N(x)$, as a function of $x=\gamma_U(u)$, for different values of     $N=N_U$.}
  \label{fig:scoring}
\end{figure}
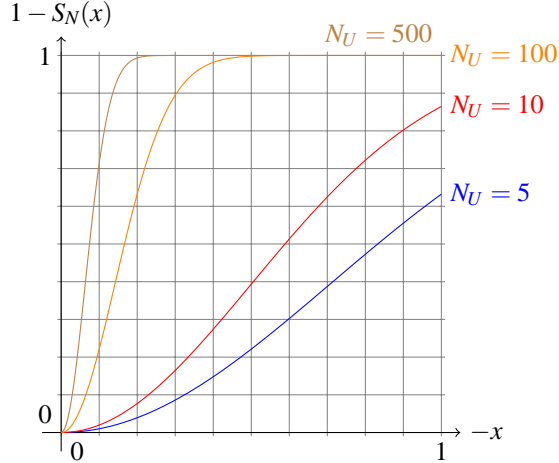

Suppose Reality follows a path up to some situation $u_o$ in $U$, which leads to an average gain $G_U(u_o)$ for Forecaster. Suppose this average gain is negative: $G_U(u_o)<0$.  We see that $\upset{u_o}\subseteq\propset{G_U<-\epsilon}$ for all $0<\epsilon<-G_U(u_o)$, and therefore all these events $\propset{G_U<-\epsilon}$ have actually occurred (because $\upset{u_o}$ has).  On the other hand, Forecaster's upper probability (in $\init$) for their occurrence satisfies $\upr(\propset{G_U<-\epsilon})\leq\exp(-\frac{N_U\epsilon^2}{4B^2})$, by Theorem~\ref{theo:largenum}. Coherence then tells us that Forecaster's upper probability (in $\init$) for the event $\upset{u_o}$, which has actually occurred, is then at most $S_{N_U}(\gamma_U(u_o))$, where
\begin{equation*}
  S_N(x)=\exp\left(-\frac{N}{4}x^2\right)
  \quad\text{and}\quad
  \gamma_U(u):=\frac{G_U(u_o)}{B}.
\end{equation*}
Observe that $\gamma_U(u_o)$ is a number in $[-1,0)$, by assumption. Coherence requires that Forecaster, because of her local predictive commitments, can be forced (by Sceptic, if he chooses his strategy well) to bet against the occurrence of the event $\upset{u_o}$ at a rate that is at least $1-S_{N_U}(\gamma_U(u_o))$. So we see that Forecaster is losing utility because of her local predictive commitments. Just how much depends on how close $\gamma_U(u_o)$ lies to $-1$ , and on how large $N_U$ is; see Fig.~\ref{fig:scoring}.
\par
\par
The upper bound $S_{N_U}(\gamma_U(u_o))$ we have constructed for the upper probability of $\upset{u_o}$ has a very interesting property, which we now try to make more explicit. Indeed, if we were to calculate Forecaster's upper probability $\upr(\upset{u_o})$ for $\upset{u_o}$ directly using Eq.~\eqref{eq:lpr}, this value would generally depend on Forecaster's predictive assessments $\rdesirs_s$ for situations $s$ that don't precede $u_o$, and that Reality therefore never got to. We shall see that such is not the case for the upper bound $S_{N_U}(\gamma_U(u_o))$ constructed using Theorem~\ref{theo:largenum}.
\par
Consider any situation $s$ before $U$ but not on the path through $u_o$, meaning that Reality never got to this situation $s$.  Therefore the corresponding gamble $h_s-m_s$ in the expression for $G_U$ isn't used in calculating the value of $G_U(u_o)$, so we can change it to anything else, and still obtain the same value of $G_U(u_o)$.
\par
Indeed, consider any other predictive model, where the only thing we ask is that the $\rdesirs'_s$ coincide with the $\rdesirs_s$ for all $s$ that precede $u_o$. For other $s$, the $\rdesirs_s'$ can be chosen arbitrarily, but still coherently. Now construct a new average gain gamble $G'_U$ for this alternative predictive model, where the only restriction is that we let $h'_s=h_s$ and $m'_s=m_s$ if $s$ precedes $u_o$. We know from the reasoning above that $G'_U(u_o)=G_U(u_o)$, so the new upper probability that the event $\upset{u_o}$ will be observed is at most
\begin{equation*}
  S_{N_U}\left(\frac{G'_U(u_o)}{B}\right)
  =S_{N_U}\left(\frac{G_U(u_o)}{B}\right)
  =S_{N_U}(\gamma_U(u_o)).
\end{equation*}
In other words, the upper bound $S_N(\gamma_U(u))$ we found for Forecaster's upper probability of Reality getting to a situation $u_o$ \emph{depends only on Forecaster's local predictive assessments $\rdesirs_s$ for situations $s$ that Reality   has actually got to, and not on her assessments for other situations}. This means that this method for scoring a predictive model satisfies Dawid's \emph{Prequential Principle}; see for instance \cite{dawid1984,dawid1999}.

\section{Concatenation and backwards recursion}\label{sec:backwards-recursion}
As we have discovered in Section~\ref{sec:concatenation}, Theorem~\ref{theo:concatenation} and Proposition~\ref{prop:local-models} enable us to calculate the global predictive lower previsions $\lpr(\cdot\vert t)$ in imprecise probability trees from local predictive lower previsions $\lpr_s$, $s\follows t$, using a backwards recursion method. That this is possible in probability trees, where the probability models are precise (previsions), is well-known,\footnote{See   Chapter~3 of Shafer's book \cite{shafer1996} on causal reasoning in probability trees. This   chapter contains a number of propositions about calculating probabilities and expectations in   probability trees that find their generalisations in Sections~\ref{sec:predictive-lower-upper}   and~\ref{sec:concatenation}.  For instance, Theorem~\ref{theo:concatenation} generalises   Proposition~3.11 in \cite{shafer1996} to imprecise probability trees.} and was arguably discovered by Christiaan Huygens in the middle of the 17-th century.\footnote{See Appendix~A of   Shafer's book \cite{shafer1996}. Shafer discusses Huygens's treatment of a special case of the   so-called \emph{Problem of Points}, where Huygens draws what is probably the first recorded   probability tree, and solves the problem by backwards calculation of expectations in the tree.   Huygens's treatment can be found in Appendix~VI of \cite{huygens16567}.} It allows for an exponential, dynamic programming-like reduction in the complexity of calculating previsions (or expectations); it seems to be essentially this phenomenon that leads to the computational efficiency of such machine learning tools as, for instance, Needleman and Wunsch's \cite{needleman1970} sequence alignment algorithm.
\par
In this section, we want to give an illustration of such exponential reduction in complexity, by looking at a problem involving Markov chains. Assume that the state $X(n)$ of a system at consecutive times $n=1,2,\dots,N$ can assume any value in a finite set $\states$. Forecaster has some beliefs about the state $X(1)$ at time $1$, leading to a coherent lower prevision $\lpr_1$ on $\gambles(\states)$. She also assesses that when the system jumps from state $X(n)=x_n$ to a new state $X(n+1)$, where the system goes to will only depend on the state $X(n)$ the system was in at time $n$, and not on the states $X(k)$ of the system at previous times $k=1,2,\dots,n-1$. Her beliefs about where the system in $X(n)=x_n$ will go to at time $n+1$ are represented by a lower prevision $\lpr_{x_n}$ on $\gambles(\states)$. 
\par
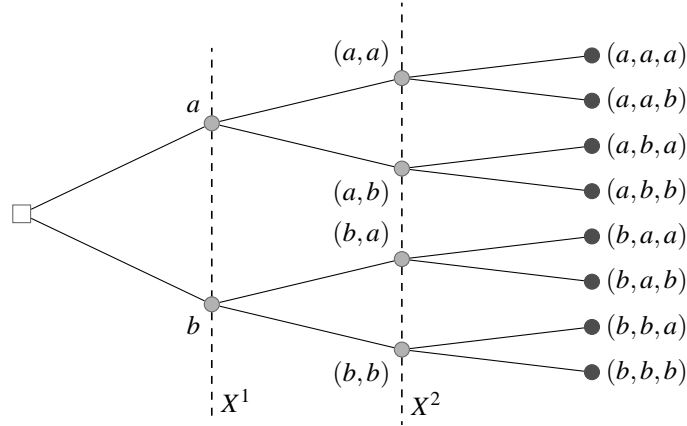
\begin{figure}[h]
  \centering
  \begin{tikzpicture}[level distance=25mm]
    \tikzstyle{root}=[rectangle,draw=black!60]   
    \tikzstyle{terminal}=[circle,draw=black!70,fill=black!70] 
    \tikzstyle{nonterminal}=[circle,draw=black!60,fill=black!30,minimum size=.5mm] 
    \tikzstyle{level 1}=[sibling distance=24mm]
    \tikzstyle{level 2}=[sibling distance=12mm]
    \tikzstyle{level 3}=[sibling distance=6mm]
    \tikzstyle{cut}=[semithick,dashed,draw=black]
    \node[root] (root) {} [grow=right,inner sep=.7mm]
    child {node[nonterminal,label=below left:{$b$}] (t) {}
      child {node[nonterminal,label=below left:{$(b,b)$}] (tt) {}
        child {node[terminal,label=right:{$(b,b,b)$}] {}}
        child {node[terminal,label=right:{$(b,b,a)$}] {}}}
      child {node[nonterminal,label=above left:{$(b,a)$}] (th) {}
        child {node[terminal,label=right:{$(b,a,b)$}] {}}
        child {node[terminal,label=right:{$(b,a,a)$}] {}}}}
    child {node[nonterminal,label=above left:{$a$}] (h) {}
      child {node[nonterminal,label=below left:{$(a,b)$}] (ht) {}
        child {node[terminal,label=right:{$(a,b,b)$}] {}}
        child {node[terminal,label=right:{$(a,b,a)$}] {}}}
      child {node[nonterminal,label=above left:{$(a,a)$}] (hh) {}
        child {node[terminal,label=right:{$(a,a,b)$}] {}}
        child {node[terminal,label=right:{$(a,a,a)$}] {}}}};
     \draw[cut] (h) -- +(0,1);
     \draw[cut] (h) -- (t) -- +(0,-1.55) node[above right] {$X^1$};
     \draw[cut] (hh) -- +(0,1);
     \draw[cut] (hh) -- (ht) -- (th) -- (tt) -- +(0,-1) node[above right] {$X^2$};
  \end{tikzpicture}
  \caption{The event tree for the time evolution of system that can be in two states, $a$ and $b$, and can change state at each time instant $n=1,2,3$. Also depicted are the respective cuts $X^1$ and $X^2$ of $\init$ where the state at times $1$ and $2$ are revealed.}
  \label{fig:markov}
\end{figure}

The time evolution of this system can be modelled as Reality traversing an event tree. An example of such a tree for $\states=\{a,b\}$ and $N=3$ is given in Fig.~\ref{fig:markov}. The situations of the tree have the form $\tuple{x}{k}\in\states^k$, $k=0,1,\dots,N$; for $k=0$ this gives some abuse of notation as we let $\states^0:=\{\init\}$. In each cut $X^k:=\states^k$ of $\init$, the value $X(k)$ of the state at time $k$ is revealed.
\par
\par
This leads to an imprecise probability tree with local predictive models $\lpr_\init:=\lpr_1$ and
\begin{equation}\label{eq:markov-condition}
  \lpr_{\tuple{x}{k}}=\lpr_{x_k}
\end{equation}
expressing the usual \emph{Markov conditional independence condition}, but here in terms of lower previsions.  For notational convenience, we now introduce a (generally non-linear) \emph{transition operator} $T$ on the linear space $\gambles(\states)$ as follows:
\begin{equation*}
  T\colon\gambles(\states)\to\gambles(\states)\colon x\mapsto\lpr_x(f),
\end{equation*}
or in other words, $T(f)$ is a gamble on $\states$ whose value $T(f)\cdot x$ in the state $x\in\states$ is given by $\lpr_x(f)$. The transition operator $T$ completely describes Forecaster's beliefs about how the system changes its state from one instant to the next.
\par
We now want to find the corresponding model for Forecaster's beliefs (in $\init$) about the state the system will be in at time $n$. So let us consider a gamble $f_n$ on $\states^N$ that actually only depends on the value $X(n)$ of $X$ at this time $n$. We then want to calculate its lower prevision $\lpr(f_n):=\lpr(f_n\vert\init)$.
\par
Consider a time instant $k\in\{0,1,\dots,n-1\}$, and a situation $\tuple{x}{k}\in\states^k$. For the children cut $C\tuple{x}{k}:=\set{(\ntuple{x}{k},x_{k+1})}{x_{k+1}\in\states}$ of $\tuple{x}{k}$, we see that $\lpr(f_n\vert C\tuple{x}{k})$ is a gamble that only depends on the value of $X(k+1)$ in $\states$, and whose value in $x_{k+1}$ is given by $\lpr(f_n\vert\ntuple{x}{k+1})$. We then find that
\begin{equation}\label{eq:markov}
  \lpr(f_n\vert\ntuple{x}{k})
  =\lpr(\lpr(f_n\vert C\tuple{x}{k})\vert\ntuple{x}{k})
  =\lpr_{x_k}(\lpr(f_n\vert C\tuple{x}{k})),
\end{equation}
where the first equality follows from Theorem~\ref{theo:concatenation}, and the second from Proposition~\ref{prop:local-models} and Eq.~\eqref{eq:markov-condition}. We first apply Eq.~\eqref{eq:markov} for $k=n-1$. By Proposition~\ref{prop:separate-coherence}.2, $\lpr(f_n\vert C\tuple{x}{n-1})=f_n$, so we are led to $\lpr(f_n\vert\ntuple{x}{n-1})=\lpr_{x_{n-1}}(f_n)=T(f_n)\cdot x_{n-1}$, and therefore
\begin{equation*}
  \lpr(f_n\vert C\tuple{x}{n-2})=T(f_n).
\end{equation*}
Substituting this in Eq.~\eqref{eq:markov} for $k=n-2$, yields $\lpr(f_n\vert\ntuple{x}{n-2})=\lpr_{x_{n-2}}(T(f_n))$, and therefore
\begin{equation*}
  \lpr(f_n\vert C\tuple{x}{n-3})=T^2(f_n).
\end{equation*}
Proceeding in this fashion until we get to $k=1$, we get $\lpr(f_n\vert C(\init))=T^{n-1}(f_n)$, and going one step further to $k=0$, Eq.~\eqref{eq:markov} yields $\lpr(f_n\vert\init)=\lpr_\init(\lpr(f_n\vert C(\init)))$ and therefore
\begin{equation}\label{eq:markov-final}
  \lpr(f_n)=\lpr_1(T^{n-1}(f_n)).
\end{equation}
\par
We see that the complexity of calculating $\lpr(f_n)$ in this way is essentially \emph{linear} in the number of time steps $n$.
\par
In the literature on imprecise probability models for Markov chains \cite{campos2003,kozine2002,skulj2006,skulj2007}, another so-called \emph{credal set}, or \emph{set of probabilities}, approach is generally used to calculate $\lpr(f_n)$. The point we want to make here is that such an approach typically has a worse (exponential) complexity in the number of time steps. To see this, recall \cite{walley1991} that a lower prevision $\lpr$ on $\gambles(\states)$ that is derived from a coherent set of really desirable gambles, corresponds to a convex closed set $\solp(\lpr)$ of probability mass functions $p$ on $\states$, called a \emph{credal set}, and given by
\begin{equation*}
  \solp(\lpr):=\set{p}{(\forall g\in\gambles(\states))\lpr(g)\leq E_p(g)}
\end{equation*}
where we let $E_p(g):=\sum_{x\in\states}p(x)g(x)$ be the expectation of the gamble $g$ associated with the mass function $p$; $E_p$ is a linear prevision in the language of Section~\ref{sec:lower-upper}. It then also holds that for all gambles $g$ on $\states$,
\begin{equation*}
  \lpr(g)
  =\min\set{E_p(g)}{p\in\solp(\lpr)}
  =\min\set{E_p(g)}{p\in\extremes\solp(\lpr)}
\end{equation*}
where $\extremes\solp(\lpr)$ is the set of extreme points of the convex closed set $\solp(\lpr)$. Typically on this approach, $\extremes(\solp(\lpr))$ is assumed to be finite, and then $\solp(\lpr)$ is called a \emph{finitely generated credal set}. See for instance \cite{cozman2000,cozman2005} for a discussion of credal sets with applications to Bayesian networks.
\par
Then $\lpr(f_n)$ can also be calculated as follows:\footnote{An explicit proof of this statement would take us to   far, but it is an immediate application of Theorems~3 and~4 in \cite{miranda2006b}.} Choose for each non-terminal situation $t=\tuple{x}{k}\in\states^k$, $k=0,1,\dots,n-1$ a mass function $p_t$ in the set $\solp(\lpr_t)$ given by Eq.~\eqref{eq:markov-condition}, or equivalently, in its set of extreme points $\extremes\solp(\lpr_t)$. This leads to a (precise) probability tree for which we can calculate the corresponding expectation of $f_n$. Then $\lpr(f_n)$ is the minimum of all such expectations, calculated for all possible assignments of mass functions to the nodes. We see that, roughly speaking, when all $\solp(\lpr_t)$ have a typical number of extreme points $M$, then the complexity of calculating $\lpr(f_n)$ will be essentially $N^n$, i.e., exponential in the number of time steps.
\par
This shows that the `lower prevision' approach can for some problems lead to more efficient algorithms than the `credal set' approach. This may be especially relevant for probabilistic inferences involving graphical models, such as credal networks \cite{cozman2000,cozman2005}.  Another nice example of this phenomenon, concerned with checking coherence for precise and imprecise probability models, is due to Walley \textit{et al.} \cite{walley2004}.

\section{Additional Remarks}
\label{sec:conclusions}
We have proved the correspondence between the two approaches only for event trees with a bounded horizon. For games with infinite horizon, the correspondence becomes less immediate, because Shafer and Vovk implicitly make use of coherence axioms that are stronger than D1--D4 and D5', leading to lower prices that dominate the corresponding predictive lower previsions. Exact matching would be restored of course, provided we could argue that these additional requirements are rational for any subject to comply with. This could be an interesting topic for further research.
\par
We haven't paid much attention to the special case that the coherent lower previsions and their conjugate upper previsions coincide, and are therefore (precise) \emph{previsions} or \emph{fair prices} in de Finetti's \cite{finetti19745} sense. When all the local predictive models $\lpr_t$ (see Proposition~\ref{prop:local-models}) happen to be precise, meaning that $\lpr_t(f)=\upr_t(f)=-\lpr_t(-f)$ for all gambles $f$ on $\wmoves_t$, then the immediate prediction model we have described in Section~\ref{sec:connections} becomes very closely related, and arguably identical to, the probability trees introduced and studied by Shafer in \cite{shafer1996}. Indeed, we then get predictive previsions $\pr(\cdot\vert s)$ that can be obtained through concatenation of the local modals $\pr_t$, as guaranteed by Theorem~\ref{theo:concatenation}.\footnote{This should for instance be compared with Proposition~3.11 in \cite{shafer1996}.}
\par
Moreover, as indicated in Section~\ref{sec:backwards-recursion}, it is possible to prove lower envelope theorems to the effect that (i) the local lower previsions $\lpr_t$ correspond to lower envelopes of sets $\solp_t$ of local previsions $\pr_t$; (ii) each possible choice of previsions $\pr_t$ in $\solp_t$ over all non-terminal situations $t$, leads to a \emph{compatible} probability tree in Shafer's \cite{shafer1996} sense, with corresponding predictive previsions $\pr(\cdot\vert s)$; and (iii) the predictive lower previsions $\lpr(\cdot\vert s)$ are the lower envelopes of the predictive previsions $\pr(\cdot\vert s)$ for the compatible probability trees. Of course, the law of large numbers of Section~\ref{sec:weak-law} remains valid for probability trees.
\par
Finally, we want to recall that Theorem~\ref{theo:concatenation} and Proposition~\ref{prop:local-models} allow for a calculation of the predictive models $\lpr(\cdot\vert s)$ using only the local models and \emph{backwards recursion}, in a manner that is strongly reminiscent of dynamic programming techniques. This should allow for a much more efficient computation of such predictive models than, say, an approach that exploits lower envelope theorems and sets of probabilities/previsions. We think that there may be lessons to be learnt from this for dealing with other types of graphical models, such as credal networks \cite{cozman2000,cozman2005}, as well.
\par
What makes this more efficient approach possible is, ultimately, the Marginal Extension Theorem (Theorem~\ref{theo:natex}), which leads to the Concatenation Formula (Theorem~\ref{theo:concatenation}), i.e., to the specific equality, rather than the general inequalities, in Proposition~\ref{prop:walley}.7. Generally speaking (see for instance \cite[Section~6.7]{walley1991} and \cite{miranda2006b}), such marginal extension results can be proved because the models that Forecaster specifies are \emph{local}, or \emph{immediate} prediction models: they relate to her beliefs, in each non-terminal situation $t$, about what move Reality is going to make \emph{immediately after} getting to $t$.  

\section*{Acknowledgements}
This paper presents research results of BOF-project 01107505. We would like to thank Enrique Miranda, Marco Zaffalon, Glenn Shafer, Vladimir Vovk and Didier Dubois for discussing and questioning some of the views expressed here, even though many of these discussions took place more than a few years ago. S\'ebastien Destercke and Erik Quaeghebeur have read and commented on earlier drafts. We are also grateful for the insightful and generous comments of three reviewers, which led us to better discuss the significance and potential applications of our results, and helped us improve the readability of this paper. 

\appendix
\section{Proofs of main results}
In this Appendix, we have gathered proofs for the most important results in the paper.
\par
We begin with a proof of Proposition~\ref{prop:walley}. Although similar results were proved for bounded gambles by Walley \cite{walley1991}, and by Williams \cite{williams2007} before him, our proof also works for the extension to possibly unbounded gambles we are considering in this paper.

\begin{proof}[Proof of Proposition~\ref{prop:walley}]
  For the first statement, we only give a proof for the first two inequalities. The proof for the remaining inequality is similar. For   the first inequality, we may assume without loss of generality that $\inf\set{\fin\in B}{f(\fin)}>-\infty$ and is therefore a real   number, which we denote by $\beta$. So we know that $I_B(f-\beta)\geq0$ and therefore $I_B(f-\beta)\in\rdesirs$, by D2. It then   follows from Eq.~\eqref{eq:lowdef} that $\beta\leq\lpr(f\vert B)$. To prove the second inequality, assume \textit{ex absurdo} that   $\upr(f\vert B)<\lpr(f\vert B)$, then it follows from Eqs.~\eqref{eq:updef} and~\eqref{eq:lowdef} that there are real $\alpha$ and   $\beta$ such that $\beta<\alpha$, $I_B(f-\alpha)\in\rdesirs$ and $I_B(\beta-f)\in\rdesirs$. By D3,   $I_B(\beta-\alpha)=I_B(f-\alpha)+I_B(\beta-f)\in\rdesirs$, but this contradicts D1, since $I_B(\beta-\alpha)<0$.
\par
We now turn to the second statement. As announced in Footnote~\ref{fn:well-defined}, we may assume that the sum of the terms $\lpr(f_1\vert B)$ and $\lpr(f_2\vert B)$ is well-defined. If either of these terms is equal to $-\infty$, the resulting inequality then holds trivially, so we may assume without loss of generality that both terms are strictly greater than $-\infty$. Consider any real $\alpha<\lpr(f_1\vert B)$ and $\beta<\lpr(f_2\vert B)$, then by Eq.~\eqref{eq:lowdef} we see that both $I_B(f_1-\alpha)\in\rdesirs$ and $I_B(f_2-\beta)\in\rdesirs$. Hence $I_B[(f_1+f_2)-(\alpha+\beta)]\in\rdesirs$, by D3, and therefore $\lpr(f_1+f_2\vert B)\geq\alpha+\beta$, using Eq.~\eqref{eq:lowdef} again.  Taking the supremum over all real $\alpha<\lpr(f_1\vert B)$ and $\beta<\lpr(f_2\vert B)$ leads to the desired inequality.
\par
To prove the third statement, first consider $\lambda>0$. Since by D4, $I_B(\lambda f-\alpha)\in\rdesirs$ if and only if $I_B(f-\alpha/\lambda)\in\rdesirs$, we get, using Eq.~\eqref{eq:lowdef}
\begin{equation*}
  \lpr(\lambda f\vert B)
  =\sup\set{\alpha}{I_B(\lambda f-\alpha)\in\rdesirs}
  =\sup\set{\lambda\beta}{I_B(f-\beta)\in\rdesirs} =\lambda\lpr(f\vert B).
\end{equation*}
For $\lambda=0$, consider that $\lpr(0\vert B)=\sup\set{\alpha}{-I_B\alpha\in\rdesirs}=0$, where the last equality follows from D1 and D2.
\par
For the fourth statement, use Eq.~\eqref{eq:lowdef} to find that
\begin{equation*}
  \lpr(f+\alpha\vert B)
  =\sup\set{\beta}{I_B(f+\alpha-\beta)\in\rdesirs}
  =\sup\set{\alpha+\gamma}{I_B(f-\gamma)\in\rdesirs} =\alpha+\lpr(f\vert B).
\end{equation*}
The fifth statement is an immediate consequence of the first.
\par
To prove the sixth statement, observe that $f_1\leq f_2$ implies that $I_B(f_2-f_1)\geq0$ and therefore $I_B(f_2-f_1)\in\rdesirs$, by D2. Now consider any real $\alpha$ such that $I_B(f_1-\alpha)\in\rdesirs$, then by D3, $I_B(f_2-\alpha)=I_B(f_1-\alpha)+I_B(f_2-f_1)\in\rdesirs$. Hence
\begin{equation*}
  \set{\alpha}{I_B(f_1-\alpha)\in\rdesirs}
  \subseteq
  \set{\alpha}{I_B(f_2-\alpha)\in\rdesirs}
\end{equation*}
and by taking suprema and considering Eq.~\eqref{eq:lowdef}, we deduce that indeed $\lpr(f_1\vert B)\leq\lpr(f_2\vert B)$.
\par
For the final statement, assume that $\lpr(f\vert C)$ is a real number for all $C\in\partit$. Also observe that $\lpr(f\vert   D)=\lpr(fI_D\vert D)$ for all non-empty $D$. Define the gamble $g$ as follows: $g(\fin):=\lpr(f\vert C)$ for all $\fin\in C$, where   $C\in\partit$. We have to prove that $\lpr(g\vert B)\leq\lpr(f\vert B)$. We may assume without loss of generality that   $\lpr(g\vert B)>-\infty$ [because otherwise the inequality holds trivially]. Fix $\epsilon>0$, and consider the gamble   $I_B(f-g+\epsilon)$.  Also consider any $C\in\partit$. If $C\subseteq B$ then $I_CI_B(f-g+\epsilon)=I_C(f-\lpr(f\vert   C)+\epsilon)\in\rdesirs$, using Eq.~\eqref{eq:lowdef}. If $C\cap B=\emptyset$ then again $I_CI_B(f-g+\epsilon)=0\in\rdesirs$, by D2.   Since $\rdesirs$ is $\partit$-conglomerable, it follows that $I_B(f-g+\epsilon)\in\rdesirs$, whence $\lpr(f-g\vert B)\geq-\epsilon$,   again using Eq.~\eqref{eq:lowdef}.  Hence $\lpr(h\vert B)\geq0$, where $h:=f-g$.  Consequently,
\begin{equation*}
  \lpr(f\vert B) 
  =\lpr(h+g\vert B) 
  \geq\lpr(h\vert B)+\lpr(g\vert B)
  \geq\lpr(g\vert B),
\end{equation*}
where we use the second statement, and the fact that $\lpr(g\vert B)>-\infty$ and $\lpr(h\vert B)\geq0$ implies that the sum on the   right-hand side of the inequality is well-defined as an extended real number.
\end{proof}

\begin{proof}[Proof of Theorem~\ref{theo:natex}]
  We have already argued that any coherent set of really desirable gambles that includes $\rdesirs$, must contain all gambles   $\gamble^\selec$ [by D3 and D5']. By D2 and D3, it must therefore include the set $\natex_\rdesirs$.  If we can show that   $\natex_\rdesirs$ is coherent, i.e., satisfies D1--D4 and D5', then we have proved that $\natex_\rdesirs$ is the natural extension   of $\rdesirs$.  This is what we now set out to do.
  \par
  We first show that D1 is satisfied. It clearly suffices to show that for no selection $\selec$, it holds that   $\gamble^\selec_\fins<0$.  This follows at once from Lemma~\ref{lem:avoiding-partial-loss} below.
  \par
  To prove that D2 holds, consider the selection $\selec_0:=0$, then $\gamble^{\selec_0}=0$, and if $f\geq0$ it follows that   $f\geq\gamble^{\selec_0}$ whence indeed $f\in\natex_\rdesirs$.
  \par
  To prove that D3 and D4 hold, consider any $f_1$ and $f_2$ in $\natex_\rdesirs$, and any non-negative real numbers $a_1$ and   $a_2$. We know there are selections $\selec_1$ and $\selec_2$ such that $f_1\geq\gamble^{\selec_1}$ and   $f_2\geq\gamble^{\selec_2}$. But $a_1\selec_1+a_2\selec_2$ is a selection as well [because the $\rdesirs_t$ satisfy D3 and D4], and   $\gamble^{a_1\selec_1+a_2\selec_2} =a_1\gamble^{\selec_1}+a_2\gamble^{\selec_2}\leq a_1f_1+a_2f_2$, whence indeed   $a_1f_1+a_2f_2\in\natex_\rdesirs$.
  \par
  To conclude, we show that D5' is satisfied. Consider any cut $U$ of $\init$.  Consider a gamble $f$ and   assume that $I_{\upset{u}}f\in\natex_\rdesirs$ for all $u\in U$. We must prove that $f\in\natex_\rdesirs$.   Let $U_{\mathit{t}}:=U\cap\fins$ and $U_{\mathit{nt}}:=U\setminus\fins$, so $U$ is the disjoint union of   $U_{\mathit{t}}$ and $U_{\mathit{nt}}$. For $\fin\in U_{\mathit{t}}$,   $I_{\upset{\fin}}f=I_{\upset{\fin}}f(\fin)\in\natex_\rdesirs$ implies that $f(\fin)\geq0$, by D1. For   $u\in U_{\mathit{nt}}$, we invoke Lemma~\ref{lem:contingency} to find that there is some $u$-selection   $\selec_u$ such that $I_{\upset{u}}f\geq\gamble^{\selec_u}_\fins$.  Now construct a selection $\selec$ as   follows.  Consider any $s$ in $\nonfins$.  If $u\precedes s$ for some [unique, because $U$ is a cut] $u\in   U_{\mathit{nt}}$, let $\selec(s):=\selec_u(s)$. Otherwise let $\selec(s):=0$.  Then
  \begin{equation*}
    \gamble^{\selec}
    =\sum_{u\in U_{\mathit{nt}}}I_{\upset{u}}\gamble^{\selec_u}
    \leq\sum_{u\in U_{\mathit{nt}}}I_{\upset{u}}f
    \leq\sum_{u\in U}I_{\upset{u}}f
    =f,
  \end{equation*}
  so indeed $f\in\natex_\rdesirs$; the first equality can be seen as immediate, or as a consequence of   Lemma~\ref{lem:decomposition}, and the second inequality holds because we have just shown that $f(\fin)\geq0$   for all $\fin\in U_t$. The rest of the proof now follows from Lemma~\ref{lem:contingency}.
\end{proof}

\begin{lemma}\label{lem:decomposition}
  Let $t$ be any non-terminal situation, and let $U$ be any cut of $t$.  Consider a $t$-selection   $\selec$, and let, for any $u\in U\setminus\fins$, $\selec_u$ be the $u$-selection given by   $\selec_u(s)=\selec(s)$ if the non-terminal situation $s$ follows $u$, and $\selec_u(s):=0$   otherwise.  Moreover, let $\selec^U$ be the $U$-called off $t$-selection for $\selec$ (as defined   after Theorem~\ref{theo:concatenation}). Then
  \begin{align*}
    \gamble^\selec_\fins &=\sum_{u\in U\cap\fins}I_{\upset{u}}\gamble^\selec(u) +\sum_{u\in U\setminus\fins}I_{\upset{u}}
    [\gamble^\selec(u)+\gamble^{\selec_u}_\fins]\\
    &=\gamble^{\selec}_U +\sum_{u\in U\setminus\fins}I_{\upset{u}}\gamble^{\selec_u}_\fins =\gamble^{\selec^U}_\fins +\sum_{u\in       U\setminus\fins}I_{\upset{u}}\gamble^{\selec_u}_\fins.
  \end{align*}
\end{lemma}

\begin{proof}
  It is immediate that the second equality holds; see Eq.~\eqref{eq:stopped-too} for the third. For the first equality, it obviously   suffices to consider the values of the left- and right-hand sides in any $\omega\in\upset{u}$ for $u\in U\setminus\fins$. The value of the right-hand side is then, using Eqs.~\eqref{eq:gamble-process} and~\eqref{eq:gamble-variable},
\begin{equation*}
  \gamble^\selec(u)+\gamble^{\selec_u}_\fins(\fin) 
  =\sum_{t\precedes s\sprecedes u}\selec(s)(u)
  +\sum_{u\precedes s\sprecedes\fin}\selec(s)(\fin)
  =\sum_{t\precedes s\sprecedes\fin}\selec(s)(\fin) 
  =\gamble^\selec_\fins(\fin).\qedhere
\end{equation*}
\end{proof}

\begin{lemma}\label{lem:avoiding-partial-loss}
  Consider any non-terminal situation $t$ and any $t$-selection $\selec$.  Then it doesn't hold that   $\gamble^\selec_\fins<0$ (on $\upset{t}$).  As a corollary, consider any cut $U$ of\/ $t$, and the gamble   $\gamble^\selec_U$ on $U$ defined by $\gamble^\selec_U(u)=\gamble^\selec(u)$. Then it doesn't hold that   $\gamble^\selec_U<0$ (on $U$).
\end{lemma}

\begin{proof}
  Define the set $P_\selec:=\set{s\in\nonfins}{\text{$t\precedes s$ and $\selec(s)\geq0$}}$, and its   (relative) complement $N_\selec:=\set{s\in\nonfins}{\text{$t\precedes s$ and $\selec(s)\not\geq0$}}$. If   $N_\selec=\emptyset$ then $\gamble^\selec_\fins\geq0$, by Eq.~\eqref{eq:gamble-variable}, so we can assume without loss of generality that   $N_\selec$ is non-empty. Consider any minimal element $t_1$ of $N_\selec$, meaning that there is no $s$   in $N_\selec$ such that $s\sprecedes t_1$ [there is such a minimal element in $N_\selec$ because of the   bounded horizon assumption]. So for all $t\precedes s\sprecedes t_1$ we have that $\selec(s)\geq0$. Choose   $\wmove_1$ in $\wmoves_{t_1}$ such that $\selec(t_1)(\wmove_1)>0$ [this is possible because   $\rdesirs_{t_1}$ satisfies D1]. This brings us to the situation $t_2:=t_1\wmove_1$. If $t_2\in N_\selec$,   then choose $\wmove_2$ in $\wmoves_{t_2}$ such that $\selec(t_2)(\wmove_2)>0$ [again possible by D1].  If   $t_2\in P_\selec$ then we know that $\selec(t_2)(\wmove_2)\geq0$ for any choice of $\wmove_2$ in   $\wmoves_{t_2}$. We can continue in this way until we reach a terminal situation   $\fin=t_1\wmove_1\wmove_2\dots$ after a finite number of steps [because of the bounded horizon   assumption].  Moreover
\begin{equation*}
  \gamble^\selec_\fins(\fin) 
  =\sum_{t\sprecedes t_1}\selec(t)(\fin(t))+\sum_{k}\selec(t_k)(\wmove_k)
  \geq0+\selec(t_1)(\wmove_1)+0>0.
\end{equation*}
It therefore can't hold that $\gamble^\selec_\fins<0$ (on $\upset{t}$).
  \par
  To prove the second statement, consider the $U$-called off $t$-selection $\selec^U$ derived from   $\selec$ by letting $\selec^U(s):=\selec(s)$ if $s$ (follows $t$ and) strictly precedes some $u$ in   $U$, and zero otherwise.  Then $\gamble^\selec(u)=\sum_{t\precedes s\sprecedes     u}\selec(s)(u)=\gamble^{\selec^U}_\fins(\fin)$ for all $\fin$ that go through $u$, where $u\in U$   [see also Eq.~\eqref{eq:stopped-too}].  Now apply the above result for the $t$-selection $\selec^U$.
\end{proof}

\begin{lemma}\label{lem:contingency}
  Consider any non-terminal situation $t$ and any gamble $f$. Then $I_{\upset{t}}f\in\natex_\rdesirs$ if and only if there is some   $t$-selection $\selec_t$ such that $I_{\upset{t}}f\geq\gamble^{\selec_t}_\fins$ (on $\upset{t}$).
\end{lemma}

\begin{proof}
  It clearly suffices to prove the necessity part. Assume therefore that $I_{\upset{t}}f\in\natex_\rdesirs$, meaning [definition of   the set $\natex_\rdesirs$] that there is some selection $\selec$ such that $I_{\upset{t}}f\geq\gamble^\selec_\fins$. Let $\selec_t$   be the $t$-selection defined by letting $\selec_t(s):=\selec(s)$ if $t\precedes s$, and zero otherwise.  It follows from   Lemma~\ref{lem:decomposition} [use the cut of $\init$ made up of $t$ and the terminal situations that do not follow $t$] that
  \begin{equation*}
    I_{\upset{t}}f
    \geq\gamble^{\selec}_\fins
    =I_{\upset{t}}[\gamble^{\selec}(t)+\gamble^{\selec_t}_\fins]
    +\sum_{\fin'\not\in\upset{t}}I_{\upset{\fin'}}\gamble^{\selec}_\fins(\fin'),
  \end{equation*}
  whence, for all $\fin\in\fins$,
  \begin{align}
    &\gamble^{\selec}_\fins(\fin)\leq0,
    \quad \fin\not\follows t\label{eq:one}\\
    &\gamble^{\selec}(t)+\gamble^{\selec_t}_\fins(\fin)\leq f(\fin), \quad \fin\follows t.\label{eq:two}
  \end{align}
  Then, by~\eqref{eq:two}, the proof is complete if we can prove that $\gamble^{\selec}(t)\geq0$.  Assume   \textit{ex absurdo} that $\gamble^{\selec}(t)<0$. Consider the cut of $\init$ made up of $t$ and the   terminal situations that don't follow $t$. Applying Lemma~\ref{lem:avoiding-partial-loss} for this cut and   for the initial situation $\init$, we see that there must be some $\fin\in\fins\setminus\upset{t}$ such   that $\gamble^{\selec}_\fins(\fin)>0$.  But this contradicts~\eqref{eq:one}.
\end{proof}

\begin{proof}[Proof of Proposition~\ref{prop:prevision-properties-general}]
  For the first statement, consider a terminal situation $\fin$ and a gamble $f$ on $\fins$. Then $\upset{\fin}=\{\fin\}$ and   therefore $I_{\upset{\fin}}(f-\alpha)=I_{\{\fin\}}(f(\fin)-\alpha)\in\natex_\rdesirs$ if and only if $\alpha\leq f(\fin)$, by D1 and   D2.  Using~Eq.~\eqref{eq:lprdef}, we find that indeed $\lpr(f\vert\fin)=f(\fin)$.  By conjugacy,   $\upr(f\vert\fin)=-\lpr(-f\vert\fin)=-(-f(\fin))=f(\fin)$ as well.
  \par
  For the second statement, use Eq.~\eqref{eq:lprdef} and observe that   $I_{\upset{t}}(f-\alpha)=I_{\upset{t}}(fI_{\upset{t}}-\alpha)$. The last statement is an immediate consequence of the second and   Proposition~\ref{prop:walley}.6.
\end{proof}

\begin{proof}[Proof of Proposition~\ref{prop:separate-coherence}]
  The first statement follows from Eq.~\eqref{eq:lprdef} if we observe that   $I_{\upset{t}}(I_{\upset{t}}-\alpha)=I_{\upset{t}}(1-\alpha)\in\natex_\rdesirs$ if and only if $\alpha\leq1$, by D1 and D2.
  \par
  For the second statement, consider any $u\in U$, then we must show that $\lpr(g\vert u)=g_U(u)$. But the $U$-measurability of $g$   tells us that $I_{\upset{u}}(g-\alpha)=I_{\upset{u}}(g_U(u)-\alpha)$, and this gamble belongs to $\natex_\rdesirs$ if and only if   $\alpha\leq g_U(u)$, by D1 and D2. Now use Eq.~\eqref{eq:lprdef}.
  \par
  The proofs of the third and fourth statements are similar, and based on the observation that   $I_{\upset{u}}(f+g-\alpha)=I_{\upset{u}}(f+g_U(u)-\alpha)$ and $I_{\upset{u}}(gf-\alpha)=I_{\upset{u}}(g_U(u)f-\alpha)$.
\end{proof}

\begin{proof}[Proof of Theorem~\ref{theo:matching}]
  First, consider an immediate prediction model $\rdesirs_t$, $t\in\nonfins$.  Define Sceptic's move spaces to be   $\smoves_t:=\rdesirs_t$ and his gain functions $\gain_t\colon\smoves_t\times\wmoves_t$ by $\gain_t(h,\wmove):=-h(\wmove)$ for all   $h\in\rdesirs_t$ and $\wmove\in\wmoves_t$. Clearly P1 and P2 are satisfied, because each $\rdesirs_t$ is a convex cone by D3   and D4. But so is the coherence requirement C. Indeed, if it weren't satisfied there would be some non-terminal situation $t$ and   some gamble $h$ in $\rdesirs_t$ such that $h(\wmove)<0$ for all $\wmove$ in $\wmoves_t$, contradicting the coherence requirement D1   for $\rdesirs_t$. We are thus led to a coherent probability protocol. We show there is matching. Consider any non-terminal situation   $t$, and any $t$-selection $\selec$. For all terminal situations $\fin\follows t$,
\begin{equation*}
  \gamble^\selec_\fins(\fin)
  =\sum_{t\precedes u\sprecedes\fin}\selec(u)(\omega(u))
  =\sum_{t\precedes u\sprecedes\fin}-\gain_u(\selec(u),\omega(u)) 
  =-\cptl^\selec_\fins(\fin),
\end{equation*}
or in other words, selections and strategies are in a one-to-one correspondence (are actually the same things), and the   corresponding gamble and capital processes are each other's inverses. It is therefore immediate from Eqs.~\eqref{eq:lprice}   and~\eqref{eq:lpr} that $\lprice_t=\lpr(\cdot\vert t)$.
\par
Conversely, consider a coherent probability protocol with move spaces $\smoves_t$ and gain functions $\gain_t\colon\smoves_t\times\wmoves_t$ for all non-terminal $t$. Define $\rdesirs_t':=\set{-\gain_t(\smove, \cdot)}{\smove\in\smoves_t}$. By a similar argument to the one above, we see that $\lpr'(\cdot\vert t)=\lprice_t$, where the $\lpr'(\cdot\vert t)$ are the predictive lower previsions associated with the sets $\rdesirs_t'$.  But each $\rdesirs_t'$ is a convex cone of gambles by P1 and P2, and by C we know that for all non-terminal situations $t$ and all gambles $h$ in $\rdesirs_t'$ there is some $\wmove$ in $\wmoves_t$ such that $h(\wmove)\geq0$. This means that the conditions for Lemma~\ref{lem:equivalence-1} are satisfied, and therefore also $\lpr'(\cdot\vert t)=\lpr(\cdot\vert t)$, where the $\lpr(\cdot\vert t)$ are the predictive lower previsions associated with the immediate prediction model $\rdesirs_t$ that is the smallest convex cone containing all non-negative gambles and including $\set{-\gain_t(\smove,\cdot)+\delta}{\smove\in\smoves_t,\delta>0}$.
\end{proof}

\begin{lemma}\label{lem:equivalence-1}
  Consider, for each non-terminal situation $t\in\nonfins$, a set of gambles $\rdesirs_t'$ on $\wmoves_t$ such that (i) $\rdesirs_t'$   is a convex cone, and (ii) for all $h\in\rdesirs_t'$ there is some $\wmove$ in $\wmoves_t$ such that $h(\wmove)\geq0$.  Then each   set $\rdesirs_t:=\set{\alpha(h+\delta)+f} {h\in\rdesirs_t',\delta>0,f\geq0,\alpha\geq0}$ is a coherent set of really desirable   gambles on $\wmoves_t$. Moreover, all predictive lower previsions obtained using the sets $\rdesirs_t$ coincide with the ones   obtained using the $\rdesirs_t'$.
\end{lemma}

\begin{proof}
  Fix a non-terminal situation $t$. We first show that $\rdesirs_t$ is a coherent set of really desirable   gambles, i.e., that D1--D4 are satisfied.  Observe that $\rdesirs_t$ is the smallest convex cone of   gambles including the set $\set{h+\delta}{h\in\rdesirs_t',\delta>0}$ and containing all non-negative   gambles. So D2--D4 are satisfied. To prove that D1 holds, consider any $g<0$ and assume \textit{ex     absurdo} that $g\in\rdesirs_t$.  Then there are $h$ in $\rdesirs_t'$, $\delta>0$, $f\geq0$ and   $\alpha\geq0$ such that $0>g=\alpha(h+\delta)+f$, whence $\alpha(h+\delta)<0$ and therefore $\alpha>0$ and   $h+\delta<0$. But by (ii), there is some $\wmove$ in $\wmoves_t$ such that $h(\wmove)\geq0$, whence   $h(\wmove)+\delta>0$. This contradicts $h+\delta<0$.
\par
We now move to the second part. Consider any gamble $f$ on $\fins$. Fix $t$ in $\nonfins$ and $\epsilon>0$. First consider any $t$-selection $\selec'$ associated with the $\rdesirs_s'$, i.e., such that $\selec'(s)\in\rdesirs_s'$ for all $s\follows t$.  Since Reality can only make a finite and \emph{bounded} number of moves, whatever happens, it is possible to choose $\delta_s>0$ for each non-terminal $s\follows t$ such that $\sum_{t\precedes s\sprecedes\fin}\delta_s<\epsilon$ for all $\fin$ in $\fins$ that follow $t$. Define the $t$-selection $\selec$ associated with the $\rdesirs_s$ by $\selec(s):=\selec'(s)+\delta_s\in\rdesirs_s$ for all non-terminal $s$ that follow $t$. Clearly $\gamble^{\selec}_\fins\leq\epsilon+\gamble^{\selec'}_\fins$, and therefore
\begin{align*}
  \lpr'(f\vert t)
  &=\sup_{\selec'}\sup\set{\alpha}{f-\alpha\geq\gamble^{\selec'}_\fins}
  \leq\sup_{\selec'}\sup\set{\alpha}{f-\alpha+\epsilon\geq\gamble^{\selec}_\fins}\\
  &=\sup_{\selec'}\sup\set{\alpha}{f-\alpha\geq\gamble^{\selec}_\fins}+\epsilon 
  \leq\lpr(f\vert t)+\epsilon.
\end{align*}
Since this inequality holds for all $\epsilon>0$, we find that $\lpr'(f\vert t)\leq\lpr(f\vert t)$.
\par
Conversely, consider any $t$-selection $\selec$ associated with the $\rdesirs_s$. For all $s\follows t$, we have that there are $h_s$ in $\rdesirs_s'$, $\delta_s>0$, $f_s\geq0$ and $\alpha_s\geq0$ such that $\selec(s)=\alpha_s(h_s+\delta_s)+f_s$. Define the $t$-selection $\selec'$ associated with the $\rdesirs_s'$ by $\selec'(s):=\alpha_sh_s=\selec(s)-\alpha_s\delta_s-f_s\leq\selec(s)$.  Clearly then also $\gamble^{\selec'}_\fins\leq\gamble^{\selec}_\fins$, and therefore
\begin{equation*}
  \lpr(f\vert t)
  =\sup_{\selec}\sup\set{\alpha}{f-\alpha\geq\gamble^\selec_\fins}
  \leq\sup_{\selec}\sup\set{\alpha}{f-\alpha\geq\gamble^{\selec'}_\fins} 
  \leq\lpr'(f\vert t).
\end{equation*}
This proves that indeed $\lpr'(f\vert t)=\lpr(f\vert t)$.
\end{proof}

\begin{proof}[Proof of Theorem~\ref{theo:concatenation}]
  It isn't difficult to see that the second statement is a consequence of the first, so we only prove the first statement.
\par
Consider any $t$-gamble $f$ on $\fins$. Recall that it is implicitly assumed that $\lpr(f\vert U)$ is again a $t$-gamble. Then we   have to prove that $\lpr(f\vert t)=\lpr(\lpr(f\vert U)\vert t)$. Let, for ease of notation, $g:=\lpr(f\vert U)$, so the $t$-gamble   $g$ is $U$-measurable, and we have to prove that $\lpr(f\vert t)=\lpr(g\vert t)$. Now, there are two possibilities.
\par
First, if $t$ is a terminal situation $\fin$, then, on the one hand, $\lpr(f\vert t)=f(\fin)$ by   Proposition~\ref{prop:prevision-properties-general}.1. On the other hand, again by   Proposition~\ref{prop:prevision-properties-general}.1,
\begin{equation*}
  \lpr(g\vert t)
  =g(\fin)
  =\lpr(f\vert U)(\fin).
\end{equation*} 
Now, since $U$ is a cut of $t=\fin$, the unique element $u$ of $U$ that $t=\fin$ goes through, is $u=\fin$, and therefore   $\lpr(f\vert U)(\fin)=\lpr(f\vert\fin)=f(\fin)$, again by Proposition~\ref{prop:prevision-properties-general}.1. This tells us that   in this case indeed $\lpr(f\vert t)=\lpr(g\vert t)$.
\par
Secondly, suppose that $t$ is not a terminal situation. Then it follows from Proposition~\ref{prop:walley}.7 and the cut conglomerability of $\natex_\rdesirs$ that $\lpr(f\vert t)\geq\lpr(\lpr(f\vert U)\vert t)=\lpr(g\vert t)$ [recall that $\lpr(\cdot\vert t)=\lpr(\cdot\vert\upset{t})$ and that $\lpr(\cdot\vert U)=\lpr(\cdot\vert\partit_U)$].  It therefore remains to prove the converse inequality $\lpr(f\vert t)\leq\lpr(g\vert t)$.  Choose $\epsilon>0$, then using Eq.~\eqref{eq:lpr} we see that there is some $t$-selection $\selec$ such that $f-\lpr(f\vert t)+\epsilon\geq\gamble^\selec_\fins$ on all paths that go through $t$.  Invoke Lemma~\ref{lem:decomposition}, using the notations introduced there, to find that
\begin{equation}\label{eq:auxiliary-concatenation-1}
  f-\lpr(f\vert t)+\epsilon
  \geq\gamble^{\selec}_U
  +\sum_{u\in U\setminus\fins}I_{\upset{u}}\gamble^{\selec_{u}}_\fins
  \quad\text{(on $\upset{t}$)}.
\end{equation}
Now consider any $u\in U$. If $u$ is a terminal situation $\fin$, then by Proposition~\ref{prop:prevision-properties-general}.1,   $g(u)=\lpr(f\vert\fin)=f(\fin)$, and therefore Eq.~\eqref{eq:auxiliary-concatenation-1} yields
\begin{equation}\label{eq:auxiliary-concatenation-2}
  g(\fin)-\lpr(f\vert t)+\epsilon
  \geq\gamble^{\selec^U}_\fins(\fin),
\end{equation}
also taking into account that $\gamble^{\selec}_U=\gamble^{\selec^U}_\fins$ [see Eq.~\eqref{eq:stopped-too}].  If $u$ is not a terminal   situation then for all $\fin\in\upset{u}$, Eq.~\eqref{eq:auxiliary-concatenation-1} yields
\begin{equation*}
  f(\fin)-\lpr(f\vert t)+\epsilon
  \geq\gamble^{\selec}_U(u) 
  +\gamble^{\selec_{u}}_\fins(\fin),
\end{equation*}
and since $\selec_u$ is a $u$-selection, this inequality together with Eq.~\eqref{eq:lpr} tells us that $\lpr(f\vert u)\geq\lpr(f\vert t)-\epsilon+\gamble^{\selec}_U(u)$, and therefore, for all $\fin\in\upset{u}$,
\begin{equation}\label{eq:auxiliary-concatenation-3}
  g(\fin)-\lpr(f\vert t)+\epsilon\geq\gamble^{\selec^U}_\fins(\fin).
\end{equation}
If we combine the inequalities~\eqref{eq:auxiliary-concatenation-2} and~\eqref{eq:auxiliary-concatenation-3}, and recall   Eq.~\eqref{eq:lpr}, we get that $\lpr(g\vert t)\geq\lpr(f\vert t)-\epsilon$.  Since this holds for all $\epsilon>0$, we may indeed   conclude that $\lpr(g\vert t)\geq\lpr(f\vert t)$.
\end{proof}

\begin{proof}[Proof of Proposition~\ref{prop:cut-reduction}]
  The condition is clearly sufficient, so let us show that it is also necessary.  Suppose that $I_{\upset{t}}f\in\natex_\rdesirs$, then   there is some $t$-selection $\selec$ such that $f\geq\gamble^\selec_\fins$, by Theorem~\ref{theo:natex} [or   Lemma~\ref{lem:contingency}]. Define, for any $u\in U\setminus\fins$, the selection $\selec_u$ as follows: $\selec_u(s):=\selec(s)$   if $s\follows u$ and $\selec_u(s):=0$ elsewhere. Then, by Lemma~\ref{lem:decomposition},
\begin{equation*}
  \gamble^\selec_\fins
  =\gamble^\selec_U+\sum_{u\in U\setminus\fins}I_{\upset{u}}\gamble^{\selec_u}_\fins.
\end{equation*}
Now fix any $u$ in $U$. If $u$ is a terminal situation $\fin$, then it follows from the equality above that
\begin{equation*}
  f_U(u)
  =f(\fin)
  \geq\gamble^\selec_U(u).
\end{equation*}
If $u$ is not a terminal situation, we get for all $\fin\in\upset{u}$:
\begin{equation*}
  f_U(u)
  =f(\fin)
  \geq\gamble^\selec_U(u)+\gamble^{\selec_u}_\fins(\fin),
\end{equation*}
whence, by taking the supremum of all $\fin\in\upset{u}$,
\begin{equation*}
  f_U(u)
  \geq\gamble^\selec_U(u)+\sup_{\fin\in\upset{u}}\gamble^{\selec_u}_\fins(\fin)
  \geq\gamble^\selec_U(u),
\end{equation*}
where the last inequality follows since $\sup_{\fin\in\upset{u}}\gamble^{\selec_u}_\fins(\fin)\geq0$ by   Lemma~\ref{lem:avoiding-partial-loss} [with $t=u$ and $\selec=\selec_u$].  Now recall that $f_U\geq\gamble^\selec_U(u)$ is   equivalent to $I_{\upset{t}}f\geq\gamble^{\selec^U}_\fins$ [see Eq.~\eqref{eq:stopped-too}].
\end{proof}

\begin{proof}[Proof of Theorem~\ref{theo:largenum}] 
  This proof builds on an intriguing idea, used by Shafer and Vovk in a different situation and form; see \cite[Lemma~3.3]{shafer2001}.
\par
Because $\abs{h_s-m_s}\leq B$ for all $t\precedes s\sprecedes u$, it follows that $G_U(u)\geq-B$, and it therefore suffices to prove the inequality for $\epsilon<B$. We work with the upper probability $\upr(\Delta_{t,\epsilon}^c\vert t)$ of the complementary event $\Delta_{t,\epsilon}^c:=\{G_U<-\epsilon\}$. It is given by
\begin{equation}\label{eq:supremum}
  \inf\set{\alpha}{\text{$\alpha-\gamble^\selec_\fins\geq
      I_{\Delta_{t,\epsilon}^c}$ for some $t$-selection $\selec$}}.
\end{equation}
Because $G_U$ is $U$-measurable, we can (and will) consider $\Delta_{t,\epsilon}^c$ as an event on $U$.  In the expression~\eqref{eq:supremum}, we may assume that $\alpha\geq0$, Indeed, if we had $\alpha<0$ and $\alpha-\gamble^\selec_\fins\geq I_{\Delta_{t,\epsilon}^c}$ for some $t$-selection $\selec$, then it would follow that $\gamble^\selec_\fins\leq\alpha<0$, contradicting Lemma~\ref{lem:avoiding-partial-loss}. Fix therefore $\alpha>0$ and $\delta>0$ and consider the selection $\selec$ such that $\selec(s):=\lambda_s(h_s-m_s)\in\rdesirs_s$ for all $t\precedes s\sprecedes U$ and let $\selec(s)$ be zero elsewhere.  Here
\begin{equation}
  \lambda_s
  :=\alpha\delta\prod_{t\precedes v\sprecedes s}[1+\delta(m_v-h_v(s))]
  =\alpha\delta\prod_{t\precedes v\sprecedes s}[1+\delta(m_v-h_v(u))],
  \label{eq:weak-law-1}
\end{equation}
where $u$ is any element of $U$ that follows $s$. Recall again that $-B\leq h_s-m_s\leq B$, so if we choose $\delta<\frac{1}{2B}$, we are certainly guaranteed that $\lambda_s>0$ and therefore indeed $\lambda_s(h_s-m_s)\in\rdesirs_s$. After some elementary manipulations we get for any $u\in U$ and any $\fin\in\upset{u}$:
\begin{equation*}
  \gamble^\selec_\fins(\fin)
  =\sum_{t\precedes s\sprecedes u}(h_s(u)-m_s)\lambda_s
  =\sum_{t\precedes s\sprecedes u}(h_s(u)-m_s)\alpha\delta 
  \prod_{t\precedes v\sprecedes     s}[1+\delta(m_v-h_v(u))]
\end{equation*}
where the second equality follows from Eq.~\eqref{eq:weak-law-1}. [The $\gamble^\selec_\fins$ is $U$-measurable.] If we let $\xi_s:=m_s-h_s(u)$ for ease of notation, then we get
\begin{align*}
  \gamble^\selec_\fins(u) 
  &=-\alpha\sum_{t\precedes s\sprecedes u}\delta\xi_s
  \prod_{t\precedes v\sprecedes s}[1+\delta\xi_v]
  =\alpha\sum_{t\precedes s\sprecedes u} 
  \prod_{t\precedes v\sprecedes s}[1+\delta\xi_v]   
  -\alpha\sum_{t\precedes s\sprecedes u}
  \prod_{t\precedes v\precedes s}[1+\delta\xi_v]\\
  &=\alpha-\alpha\prod_{t\precedes v\sprecedes u}[1+\delta\xi_v]
  =\alpha-\alpha\prod_{t\precedes v\sprecedes u}[1+\delta(m_v-h_v(u))]
\end{align*}
for all $u$ in $U$.  Then it follows from~\eqref{eq:supremum} that if we can find an $\alpha\geq0$ such that
\begin{equation*}
  \alpha\prod_{t\precedes v\sprecedes u}[1+\delta(m_v-h_v(u)))]\geq1
\end{equation*}
whenever $u$ belongs to $\Delta_{t,\epsilon}^c$, then this $\alpha$ is an upper bound for $\upr(\Delta_{t,\epsilon}^c\vert t)$.  By   taking logarithms on both sides of the inequality above, we get the equivalent condition
\begin{equation}\label{eq:equivalent}
  \ln\alpha+\sum_{t\precedes s\sprecedes u}\ln[1+\delta(m_s-h_s(u))]\geq0.
\end{equation}
Since $\ln(1+x)\geq x-x^2$ for $x>-\frac{1}{2}$, and $\delta(m_s-h_s(u))\geq-\delta B>-\frac{1}{2}$ by our previous restrictions on $\delta$, we find
\begin{align*}
  \sum_{t\precedes s\sprecedes u}\ln[1+\delta(m_s-h_s(u))] &\geq\sum_{t\precedes s\sprecedes     u}\delta(m_s-h_s(u))
-\sum_{t\precedes s\sprecedes u}[\delta(m_s-h_s(u))]^2\\
&\geq\delta\sum_{t\precedes s\sprecedes u}
[m_s-h_s(u)]-\delta^2n_U(u)B^2\\
&=n_U(u)\delta\left[-G_U(u)-B^2\delta\right].
\end{align*}
But for all $u\in\Delta_{t,\epsilon}^c$, $-G_U(u)>\epsilon$, so for all such $u$
\begin{equation*}
  \sum_{t\precedes s\sprecedes u}\ln[1+\delta(m_s-h_s(u))]
  >n_U(u)\delta(\epsilon-B^2\delta).
\end{equation*}
If we therefore choose $\alpha$ such that for all $u\in U$, $\ln\alpha+n_U(u)\delta(\epsilon-B^2\delta)\geq0$, or equivalently   $\alpha\geq\exp(-n_U(u)\delta(\epsilon-B^2\delta))$, then the above condition~\eqref{eq:equivalent} will indeed be satisfied for all   $u\in\Delta_{t,\epsilon}^c$, and then $\alpha$ is an upper bound for $\upr(\Delta_{t,\epsilon}^c\vert t)$. The tightest (smallest) upper bound is always (for all $u\in U$) achieved for $\delta=\frac{\epsilon}{2B^2}$. Replacing $n_U$ by   its minimum $N_U$ allows us to get rid of the $u$-dependence, so we see that   $\upr(\Delta_{t,\epsilon}^c\vert t)\leq\exp(-\frac{N_U\epsilon^2}{4B^2})$.  We previously required that $\delta<\frac{1}{2B}$, so if we use   this value for $\delta$, we find that we have indeed proved this inequality for $\epsilon<B$.
\end{proof}


\begin{thebibliography}{10} 
\bibitem{boole1847} G.~Boole. \newblock {\em The Laws of Thought}. \newblock Dover Publications, New York, 1847, reprint 1961. 
\bibitem{campos2003} M.~A. Campos, G.~P. Dimuro, A.~C. {d}a Rocha~Costa, and V.~Kreinovich. \newblock Computing 2-step predictions for interval-valued finite stationary   {M}arkov chains. \newblock Technical Report UTEP-CS-03-20a, University of Texas at El Paso,   2003. 
\bibitem{cozman2000} F.~G. Cozman. \newblock Credal networks. \newblock {\em Artificial Intelligence}, 120:199--233, 2000. 
\bibitem{cozman2005} F.~G. Cozman. \newblock Graphical models for imprecise probabilities. \newblock {\em International Journal of Approximate Reasoning},   39(2-3):167--184, June 2005. 
\bibitem{dawid1984} A.~Ph. Dawid. \newblock Statistical theory: The prequential approach. \newblock {\em Journal of the Royal Statistical Society, Series A},   147:278--292, 1984. 
\bibitem{dawid1999} A.~Ph. Dawid and V.~G. Vovk. \newblock Prequential probability: principles and properties. \newblock {\em Bernoulli}, 5:125--162, 1999. 
\bibitem{cooman2007b} G.~{d}e Cooman and F.~Hermans. \newblock On coherent immediate prediction: Connecting two theories of   imprecise probability. \newblock In G.~{d}e Cooman, J.~Vejnarova, and M.~Zaffalon, editors, {\em   ISIPTA '07 -- Proceedings of the Fifth International Symposium on Imprecise   Probability: Theories and Applications}, pages 107--116. SIPTA, 2007.
\bibitem{cooman2005c} G.~{d}e Cooman and E.~Miranda. \newblock Symmetry of models versus models of symmetry. \newblock In W.~L. Harper and G.~R. Wheeler, editors, {\em Probability and   Inference: Essays in Honor of {Henry E.~Kyburg, Jr.}}, pages 67--149. King's   College Publications, 2007.
\bibitem{cooman2004b} G.~{d}e Cooman and M.~Zaffalon. \newblock Updating beliefs with incomplete observations. \newblock {\em Artificial Intelligence}, 159(1-2):75--125, November 2004. 
\bibitem{finetti1970} B.~{d}e Finetti. \newblock {\em Teoria delle Probabilit\`a}. \newblock Einaudi, Turin, 1970. 
\bibitem{finetti19745} B.~{d}e Finetti. \newblock {\em Theory of Probability: A Critical Introductory Treatment}. \newblock John Wiley \& Sons, Chichester, 1974--1975. \newblock {E}nglish translation of \cite{finetti1970}, two volumes. \bibitem{gardenfors1988a} P.~G\"ardenfors and N.-E. Sahlin. \newblock {\em Decision, Probability, and Utility}. \newblock Cambridge University Press, Cambridge, 1988. 
\bibitem{goldstein1983} M.~Goldstein. \newblock The prevision of a prevision. \newblock {\em Journal of the American Statistical Society}, 87:817--819, 1983. \bibitem{hoeffding1963} W.~Hoeffding. \newblock Probability inequalities for sums of bounded random variables. \newblock {\em Journal of the Americal Statistical Association}, 58:13--30,   1963. 
\bibitem{huygens16567} Ch. Huygens. \newblock {\em Van Rekeningh in Spelen van Geluck}. \newblock 1656--1657. \newblock Reprinted in Volume XIV of \cite{huygens1888}. 
\bibitem{huygens1888} Ch. Huygens. \newblock {\em {\OE}uvres compl\`etes de Christiaan Huygens}. \newblock Martinus Nijhoff, Den Haag, 1888-1950. \newblock Twenty-two volumes. Available in digitised form from the   Biblioth\`eque nationale de France (\texttt{http://gallica.bnf.fr}). 
\bibitem{kozine2002} Igor~O. Kozine and Lev~V. Utkin. \newblock Interval-valued finite markov chains. \newblock {\em Reliable Computing}, 8(2):97--113, April 2002. 
\bibitem{kyburg1964} H.~E. Kyburg~Jr. and H.~E. Smokler, editors. \newblock {\em Studies in Subjective Probability}. \newblock Wiley, New York, 1964. \newblock Second edition (with new material) 1980. 
\bibitem{manski2003} C.~Manski. \newblock {\em Partial Identification of Probability Distributions}. \newblock Springer-Verlag, New York, 2003.
\bibitem{miranda2006b} E.~Miranda and G.~{d}e Cooman. \newblock Marginal extension in the theory of coherent lower previsions. \newblock {\em International Journal of Approximate Reasoning}, 46(1):188--225,   September 2007. 
\bibitem{needleman1970} S.~B. Needleman and C.~D. Wunsch. \newblock A general method applicable to the search for similarities in the   amino acid sequence of two proteins. \newblock {\em Journal of Molecular Biology}, 48:443--453, 1970. 
\bibitem{ramsey1931} F.~P. Ramsey. \newblock Truth and probability (1926). \newblock In R.~B. Braithwaite, editor, {\em The Foundations of Mathematics and   other Logical Essays}, chapter VII, pages 156--198. Kegan, Paul, Trench,   Trubner \& Co., London, 1931. \newblock Reprinted in \cite{kyburg1964} and \cite{gardenfors1988a}. 
\bibitem{shafer1982} G.~Shafer. \newblock Bayes's two arguments for the {R}ule of {C}onditioning. \newblock {\em The Annals of Statistics}, 10:1075--1089, 1982. 
\bibitem{shafer1983} G.~Shafer. \newblock A subjective interpretation of conditional probability. \newblock {\em Journal of Philosophical Logic}, 12:453--466, 1983. 
\bibitem{shafer1985} G.~Shafer. \newblock Conditional probability. \newblock {\em International Statistical Review}, 53:261--277, 1985.
\bibitem{shafer1996a} G.~Shafer. \newblock {\em The Art of Causal Conjecture}. \newblock The MIT Press, Cambridge, MA, 1996. 
\bibitem{shafer1996} G.~Shafer. \newblock The significance of {J}acob {B}ernoulli's \emph{Ars Conjectandi} for   the philosophy of probability today. \newblock {\em Journal of Econometrics}, 75:15--32, 1996. 
\bibitem{shafer2000} G.~Shafer, P.~R. Gillett, and R.~Scherl. \newblock The logic of events. \newblock {\em Annals of Mathematics and Artificial Intelligence}, 28:315--389,   2000. 
\bibitem{shafer2003} G.~Shafer, P.~R. Gillett, and R.~B. Scherl. \newblock A new understanding of subjective probability and its generalization   to lower and upper prevision. \newblock {\em International Journal of Approximate Reasoning}, 33:1--49, 2003. 
\bibitem{shafer2001} G.~Shafer and V.~Vovk. \newblock {\em Probability and Finance: It's Only a Game!} \newblock Wiley, New York, 2001. \bibitem{vovk2005} V.~Vovk, A.~Gammerman, and G.~Shafer. \newblock {\em Algorithmic learning in a Random World}. \newblock Springer, New York, 2005. 
\bibitem{skulj2006} D.~\v{S}kulj. \newblock Finite discrete time {M}arkov chains with interval probabilities. \newblock In J.~Lawry, E.~Miranda, A.~Bugarin, S.~Li, M.~A. Gil,   P.~Grzegorzewski, and O.~Hryniewicz, editors, {\em Soft Methods for   Integrated Uncertainty Modelling}, pages 299--306. Springer, Berlin, 2006. 
\bibitem{skulj2007} D.~\v{S}kulj. \newblock Regular finite {M}arkov chains with interval probabilities. \newblock In G.~{d}e Cooman, J.~Vejnarova, and M.~Zaffalon, editors, {\em   ISIPTA '07 -- Proceedings of the Fifth International Symposium on Imprecise   Probability: Theories and Applications}, pages 405--413. SIPTA, 2007. 
\bibitem{walley1991} P.~Walley. \newblock {\em Statistical Reasoning with Imprecise Probabilities}. \newblock Chapman and Hall, London, 1991. \bibitem{walley1996} P.~Walley. \newblock Measures of uncertainty in expert systems. \newblock {\em Artificial Intelligence}, 83(1):1--58, May 1996.
\bibitem{walley2000} P.~Walley. \newblock Towards a unified theory of imprecise probability. \newblock {\em International Journal of Approximate Reasoning}, 24:125--148,   2000. 
\bibitem{walley2004} P.~Walley, R.~Pelessoni, and P.~Vicig. \newblock Direct algorithms for checking consistency and making inferences from   conditional probability assessments. \newblock {\em Journal of Statistical Planning and Inference}, 126:119--151,   2004. 
\bibitem{wasserman2004} L.~Wasserman. \newblock {\em All of Statistics}. \newblock Springer, New York, 2004. 
\bibitem{williams1975} P.~M. Williams. \newblock Notes on conditional previsions. \newblock Technical report, School of Mathematical and Physical Science,   University of Sussex, UK, 1975. 
\bibitem{williams2007} P.~M. Williams. \newblock Notes on conditional previsions. \newblock {\em International Journal of Approximate Reasoning}, 44:366--383,   2007. \newblock Revised journal version of \cite{williams1975}. \end{thebibliography}

\end{document}